\documentclass[12pt]{amsart}
\usepackage[top=1.08in, left=1.28in, bottom=1.08in, right=1.28in]{geometry}
\usepackage{amsmath,amsfonts,amsbsy,amsgen,amscd,mathrsfs,amssymb,amsthm,mathtools,bbm,enumerate}
\usepackage[colorlinks=true,citecolor=teal!80!black,linkcolor=blue]{hyperref}
\usepackage[foot]{amsaddr}

\usepackage{hhline}
\usepackage{fouriernc}

\usepackage{pgfplots}
\usepackage{tikz-cd}
\usetikzlibrary{arrows,automata}

\numberwithin{equation}{section}
\newtheorem{theorem}{Theorem}[section]
\newtheorem{corollary}[theorem]{Corollary}
\newtheorem{conjecture}[theorem]{Conjecture}
\newtheorem{lemma}[theorem]{Lemma}
\newtheorem{proposition}[theorem]{Proposition}
\theoremstyle{definition}
\newtheorem{assumption}[theorem]{Assumption}
\newtheorem{definition}[theorem]{Definition}

\newtheorem{notation}[theorem]{Notation}
\newtheorem{problem}[theorem]{Problem}
\newtheorem{openproblem}[theorem]{Open Problem}

\newtheorem{remark}[theorem]{Remark}

\makeatletter

\newcommand\al{\alpha}
\newcommand\be{\beta}
\newcommand\dd{\mathrm d}
\newcommand\De{\Delta}
\newcommand\de{\delta}
\newcommand\deq{\stackrel{\mathrm{distr.}}{=}}
\newcommand\eps{\varepsilon}
\newcommand\ga{\gamma}
\newcommand\Ga{\Gamma}
\newcommand\ka{\kappa}

\newcommand{\om}{\omega}

\newcommand\Si{\Sigma}
\newcommand\si{\sigma}

\newcommand\ze{\zeta}

\renewcommand\bar{\overline}
\renewcommand\d{~\mathrm d}
\renewcommand\phi{\varphi}
\renewcommand\rho{\varrho}
\renewcommand\th{\vartheta}
\renewcommand\hat{\widehat}

\newcommand\bs{\boldsymbol}
\newcommand\mbb{\mathbb}
\newcommand\mbf{\mathbf}
\newcommand\mc{\mathcal}
\newcommand\mf{\mathfrak}
\newcommand\mr{\mathrm}

\newcommand\msf{\mathsf}

\begin{document}

\title[The PAM's Total Mass at Small Times]{The Parabolic Anderson Model's Total Mass at Small Times: Geometry, Fluctuations, and Renormalization}
\author{Pierre Yves Gaudreau Lamarre}
\address{Department of Mathematics,
Syracuse University,
Syracuse, NY 13244}
\email{pgaudrea@syr.edu}
\author{Yuanyuan Pan}
\email{Yuanypan@hotmail.com}
\date{\today}
\maketitle

\begin{abstract}
Let $D\subset\mathbb R^d$ be a bounded domain,
let $\kappa>0$ be fixed, and let $W$ be a fractional Brownian sheet on $\mathbb R\times\mathbb R^d$.
Consider the Stratonovich parabolic Anderson model (PAM) $\partial_tu_\kappa=(\frac12\Delta+\kappa W')u_\kappa$
with Dirichlet boundary condition on $D$ and the flat initial condition $u_\kappa(0,\cdot)=\mathbf 1_D$.
We calculate exact asymptotics for the expectation and the standard deviation of the total mass $\int_Du_\kappa(t,x)~\mathrm d x$
as $t\to0$ under the assumption that
$W$'s Hurst indices are all at least $1/2$ and that $u_\kappa$'s moments are finite for small enough $t>0$.
In doing so, we uncover that these asymptotics are determined by a competition between three mechanisms:\\
(1) {\bf Geometry:} The rate of heat diffusion through the boundary $\partial D$.\\
(2) {\bf Fluctuations:} $W$'s time Hurst index.\\
(3) {\bf Renormalization:} The singularity of deterministic Stratonovich corrections.\\
As a result, we identify novel phase transition phenomena, which arise from
the influence of $W$'s Hurst indices on the relative magnitudes of these contributions.
\end{abstract}

\section{Introduction}

\subsection{Setup}

Let $d\in\mbb N$, and let $W$ be a fractional Brownian sheet on
$\mbb R\times\mbb R^d$ with Hurst parameter $H_0$ in time
and Hurst parameters $H_1,\ldots,H_d$ in space.
Then, let
\[\xi(t,x)=\frac{\partial^{d+1}}{\partial t\partial x_1\cdots\partial x_d}W(t,x),\qquad (t,x)\in\mbb R\times\mbb R^d\]
be the distributional derivative of $W$. If $H_0=1$, then $\xi$ is time-independent.
Let $D\subset\mbb R^d$ be a bounded domain (i.e., nonempty, open, and connected).
For every $\ka\geq0$, let $u_\ka$ be the
Stratonovich parabolic Anderson model (PAM) with Dirichlet boundary condition on $D$ and noise $\ka\xi$.
That is, the formal solution of the stochastic PDE
\begin{align}
\label{Equation: Dirichlet PAM}
\partial_tu_\ka(t,x)=\big(\tfrac12\De+\ka\xi(t,x)\big)u_\ka(t,x),\qquad t>0\text{ and }x\in D,
\end{align}
where $u_\ka(t,x)=0$ for every $t>0$ and $x\in\partial D$.
Throughout the paper,
we assume the flat initial condition $u_\ka(0,x)=1$ for $x\in D.$

When $\ka=0$, \eqref{Equation: Dirichlet PAM} reduces to the classical Dirichlet heat equation on $D$.
Broadly speaking, we are interested in understanding the effect that the
transition from $\ka=0$ to $\ka>0$ (i.e., the introduction of noise)
has on the PAM's evolution in time. 
At the time of writing this paper,
the PAM literature is mostly
concerned with the effect of $\ka$ on large-$t$ asymptotics. The main motivation for this
is to study the intermittency phenomenon;
see, e.g., \cite{ChenHuKalbasiNualart,ChenQuenched,
ChenFractionalSkorokhod,ChenRoughSpace,ChenRoughTime,
ChenHuSongSong,HuHuangNualartTindel}.
In sharp contrast to this, in this paper we are concerned with
the following:

\begin{problem}
\label{Problem: Main}
Given a fixed $\ka>0$, calculate the small-$t$ asymptotics of the total mass
perturbation function
\begin{align}
\label{Equation: Qkappa}
\msf Q_\ka(t)=\int_{D}u_\ka(t,x )\d x-\int_{D}u_0(t,x )\d x,\qquad t>0.
\end{align}
In particular, understand how different choices of the parameters
$H_0,H_1,\ldots,H_d$ influence this small-time behavior.
\end{problem}

\subsection{Main Result}

In \cite{GaudreauLamarrePan}, we studied Problem \ref{Problem: Main}
in the special case where $\xi$ is a time-independent
white noise on $\mbb R^2$ (i.e., $d=2$, $H_0=1$, and $H_1=H_2=1/2$).
In \cite[Theorem 1.3]{GaudreauLamarrePan}, we proved an exact asymptotic for
$\mbf E\big[\msf Q_\ka(t)\big]$ and a non-explicit upper bound for
$\mbf{Var}\big[\msf Q_\ka(t)\big]$, both as $t\to0$.
See also \cite[Section 2.2]{GaudreauLamarrePan} for applications
of this result in spectral geometry.
In this paper, we provide exact asymptotics for
both $\mbf E\big[\msf Q_\ka(t)\big]$ and $\sqrt{\mbf{Var}[\msf Q_\ka(t)]}$
for a much wider class of fractional noises:

\begin{assumption}
\label{Assumption}
Let $1/2\leq H_0\leq 1$ and
$1/2\leq H_1,\ldots,H_d<1$. Moreover,
if we denote $H=H_1+\cdots+H_d$, then we assume that
\begin{align}
\label{Equation: Assumption}
d-H\leq1\qquad\text{and}\qquad d-H<2H_0-1/2.
\end{align}
\end{assumption}

We make this assumption for two reasons: Firstly, since $H_i\geq1/2$ for all $i\geq0$,
we can formally view $\xi$ as a centered Gaussian process with nonnegative covariance
\[\mbf E\big[\xi(s,x)\xi(t,y)\big]=\msf h(s-t)\msf g(x-y),\qquad (s,x),(t,y)\in\mbb R\times \mbb R^d,\]
where, if we let $J=\{1\leq j\leq d:H_j=1/2\}$ and $\de_0$ be the delta Dirac distribution, then
\begin{align}
\label{Equation: g}
\msf g(x-y)=\left(\prod_{j\in J}\de_0(x_j-y_j)\right)\left(\prod_{j\not\in J}H_j(2H_j-1)|x_j-y_j|^{2H_j-2}\right),
\qquad x,y\in\mbb R^d,
\end{align}
and
\begin{align}
\label{Equation: h}
\msf h(s-t)=\begin{cases}
H_0(2H_0-1)|s-t|^{2H_0-2}&H_0>1/2,\\
\de_0(s-t)&H_0=1/2,
\end{cases}
\qquad s,t\in\mbb R.
\end{align}
Secondly, and most importantly, \eqref{Equation: Assumption}
ensures that the total mass' moments are finite for small enough $t>0$.
See Proposition \ref{Proposition: L^2} and Definition \ref{Definition: L^2} for
a simple construction of the total mass (for a fixed $t$) as an $L^2$ limit under these assumptions.

In this context, our main result is as follows:

\begin{notation}
Given two nonvanishing functions $f$ and $g$ defined in a neighborhood of zero, we use
"$f(t)\sim g(t)$" as a shorthand for "$f(t)/g(t)\to1$ as $t\to0$."
\end{notation}

\begin{notation}
Let $|D|$ denote $D$'s volume.
Define the constants
\begin{align}
\label{Equation: rho exponent}
&\rho=(2H_0-1)-(d-H),\\
\label{Equation: b Constant}
&\mf b=\frac{H_0(2H_0-1)}{2^{H}}\prod_{j=1}^d\frac{\Ga(2H_j+1)}{\Ga(H_j)},\\
\label{Equation: c Constant}
&\mf c=\left(\int_{D^2}\msf g(x-y)\d x\dd y\right)^{1/2}.
\end{align}
(Note that, under Assumption \ref{Assumption}, $\rho>-1/2$ because $d-H<2H_0-1/2$, and $\rho<1$ because $d-H>0$ and $2H_0-1\leq1$. Moreover, $\mf b\geq0$ because $H_0\geq1/2$.)
\end{notation}

\begin{theorem}
\label{Theorem: Main}
Suppose that Assumption \ref{Assumption} holds, and let $\ka>0$.
On the one hand:
\begin{enumerate}
\item If $0<\rho<1$, then
$\displaystyle\mbf E\big[\msf Q_\ka(t)\big]\sim\frac{\ka^2\mf b|D|}{(1+\rho)\rho}\,t^{1+\rho}.$
\vspace{5pt}
\item If $\rho=0$, then
$\displaystyle\mbf E\big[\msf Q_\ka(t)\big]\sim-\ka^2\mf b|D|\,t\log(1/t).$
\vspace{5pt}
\item If $-1/2<\rho<0$, then
$\displaystyle\mbf E\big[\msf Q_\ka(t)\big]
\begin{cases}
\sim\frac{\ka^2\mf b|D|}{(1+\rho)\rho}\,t^{1+\rho}
&\text{if }H_0>1/2,\\
=0\quad\text{for all }t>0&\text{if }H_0=1/2.
\end{cases}$
\end{enumerate}
On the other hand, in all cases,
\begin{align}
\label{Equation: Variance}
\sqrt{\mbf{Var}[\msf Q_{\ka}(t)]}\sim \ka\mf c\,t^{H_0}.
\end{align}
\end{theorem}

\begin{remark}
\label{Remark: b vanishing and rho singularity}
The $H_0=1/2$ exception in Theorem \ref{Theorem: Main}-(3)
can be explained by the fact that the constant $\mf b$ in \eqref{Equation: b Constant}, which appears in every expectation asymptotic, vanishes in that case (due to
the presence of the factor $2H_0-1$).
Moreover, the appearance of the logarithm $\log(1/t)$ when $\rho=0$ coincides with a singularity of the constant $\frac{1}{(1+\rho)\rho}$, the latter of which also changes sign when $\rho$ goes from positive to negative.
\end{remark}

\subsection{Interpretation}

The study of the small-time asymptotics of the quantity
\[\int_{D}u_0(t,x)\d x,\] known
as the heat content of $D$, has a distinguished history
in analysis and geometry. See, e.g.,
\cite{Gilkey3,Gilkey2}, and
\cite[Chapter 2]{Gilkey}
for thorough surveys.
In short, the main goal of this theory is to understand how various features
of $D$'s geometry influence the rate at which heat dissipates through the boundary
$\partial D$.

The effect of perturbations by regular environments on the heat content
 always has the same general form. For instance, if $\xi$ were continuous, then
a straightforward calculation (e.g., using Feynman-Kac) would imply that
\begin{align}
\label{Equation: Smooth Case}
\int_Du_\ka(t,x)\d x=\int_D u_0(t,x)\d x+\ka\left(\int_D\xi(0,x)\d x\right)t+o(t)\qquad\text{as }t\to0.
\end{align}
See also \cite[Chapter 2]{Gilkey} for a complete power series expansion when $\xi$ is smooth.
It is known that the heat content cannot converge to $|D|$ (which is its $t\to0$ limit) at a faster rate than $t^{1/2}$; see \cite{BurchardSchmuckenschlager}. In particular,
regular environments can never compete with the leading scale in \eqref{Equation: Smooth Case}.
Thus, our motivation for Problem \ref{Problem: Main} comes entirely
from the potential for more interesting asymptotics when $\xi$ is singular.

In \cite{GaudreauLamarrePan}, we observed one isolated instance of this:
In addition to a random term of order at most $t$,
the two-dimensional white noise PAM
has a deterministic contribution of order $t\log(1/t)$ in its total mass.
This then leads us to the main insight of this paper:
Thanks to the two improvements on \cite{GaudreauLamarrePan} in Theorem \ref{Theorem: Main}
(i.e., the generality of Assumption \ref{Assumption}
and the optimality of the variance asymptotic in \eqref{Equation: Variance}),
our main result sheds new light on the true diversity of asymptotics that are
possible with singular noises. Moreover, our results identify the mechanisms responsible
for different types of asymptotics.

More specifically, if $\xi=W'$ for a fractional Brownian sheet $W$, then an informal Stratonovich chaos expansion
of the total mass suggests that
\begin{align}
\label{Equation: Formal Chaos 0}
\int_Du_\ka(t,x)\d x=\int_Du_0(t,x)\d x+\sum_{n=1}^\infty \ka^nI^\circ_n(t),
\end{align}
where $I_n^\circ(t)$ is the iterated Stratonovich stochastic integral
that appears in the $n^{\mr{th}}$ iteration of Duhamel's principle
(see, e.g., \cite[Section 2]{Bal} and \cite[Section 4.1]{Bal2}). By the Hu-Meyer conversion formula
(\cite[(5)]{HuMeyer}),
we expect that we can write
\[I^\circ_n(t)=\sum_{m=0}^{\lfloor n/2\rfloor}I^{(n)}_{n-2m}(t)\qquad\text{for every }n\geq1,\]
where $I^{(n)}_{n-2m}(t)$ is a $(n-2m)$-fold iterated Wiener-It\^o integral (with some combinatorial coefficient).
In particular, if we isolate the zero-order contributions in these sums
(which only occur when $n$ is even),
then we get the heuristic decomposition
\begin{align}
\label{Equation: Formal Chaos}
\int_Du_\ka(t,x)\d x=\underbrace{\int_Du_0(t,x)\d x}_{\text{heat content}}
+\underbrace{\sum_{n=1}^\infty \ka^n\msf F_n(t)}_{\text{fluctuations}}
+\underbrace{\sum_{n=1}^\infty \ka^{2n}\msf D_{2n}(t)}_{\text{deterministic}},
\end{align}
where, for all $n\geq1$, we define
\[\msf F_n(t)=\begin{cases}
I^\circ_n(t)&n\text{ odd},\\
I^\circ_n(t)-I^{(n)}_0(t)&n\text{ even},
\end{cases}
\qquad\text{and}\qquad
\msf D_{2n}(t)=I^{(2n)}_0(t).
\]

Thus, the total mass asymptotics are determined by a competition between
the three contributions on the right-hand side of \eqref{Equation: Formal Chaos}.
The first is well-understood:
\begin{enumerate}
\item {\bf Geometry:} The rate at which
the heat content converges to $|D|$ is determined by the geometric properties of $\partial D$.
This is typically on the order of $t^q$ for some $0<q\leq 1/2$. See Section \ref{Section: Geometry} and references therein for details.
\end{enumerate}
Beyond this classical theory, our main result identifies the mechanisms that control
the other two contributions when Assumption \ref{Assumption} holds:
\begin{enumerate}
 \setcounter{enumi}{1}
 \item {\bf Fluctuations:} Since the space variable is integrated in the total mass,
the space Hurst indices $H_1,\ldots,H_d$ only contribute to the fluctuations
via the constant $\mf c$ in \eqref{Equation: c Constant}. In particular, the size
of the fluctuations term is entirely determined by $t^{H_0}$, which is the scale of
$W$'s standard deviation in time.
\item {\bf Renormalization:} The size of the deterministic contribution is controlled by the
constant $\rho$ in \eqref{Equation: rho exponent}.
As per Theorem \ref{Theorem: Main}, the deterministic
term can only compete with the heat content or
the fluctuations when $\rho\leq0$.
The parameter $\rho$ also characterizes whether or not the PAM needs a renormalization, the latter of which is
necessary if and only if
$\rho\leq0$ under Assumption \ref{Assumption} (see Definition \ref{Definition: Renormalization} and Proposition \ref{Proposition: L^2}).
In fact, in Theorem \ref{Theorem: Main}-(2) and -(3) (except when $H_0=1/2$),
the decay rate of the expectation matches exactly the renormalization
singularity, up to an additional factor of $t$. See
Corollary \ref{Corollary: Renormalization} and Remark \ref{Remark: eps/t 1} for a more detailed explanation of this phenomenon.
\end{enumerate}
One of the most interesting consequences of this insight is that we can identify
two critical thresholds at which any change in $W$'s parameters leads to a total mass "phase transition"
(in the sense of a change in the nature of the dominant small-time mechanism in \eqref{Equation: Formal Chaos}),
namely:
\begin{enumerate}
\item If $H_0=1/2$, then the fluctuations term can catch up to the heat content.
\item If $\rho=H_0-1$, then the deterministic and fluctuations terms are the same size.
\end{enumerate}
We point to Corollaries \ref{Corollary: Geometry-Fluctuations} and \ref{Corollary: Fluctuations-Renormalization}
for detailed statements,
as well as Figure \ref{Figure: Phase Diagram} for a complete phase diagram
involving the relevant parameters.

\subsection{Beyond Assumption \ref{Assumption}}
\label{Section: Limitations}

In recent years, the fractional PAM has been studied under less restrictive
conditions than Assumption \ref{Assumption}.
Notably, there has been a recent interest in the setting of so-called
rough fractional noise, which corresponds to $H_i<1/2$;
see, e.g., \cite{ChenRoughSpace,ChenRoughTime,ChenDeyaOuyangTindel2}.
Moreover, the PAM has been constructed in cases where the requirement \eqref{Equation: Assumption}
is not satisfied; see \cite{HairerKPZ,HairerLabbe}. It is natural to ask what
becomes of the three contributions in \eqref{Equation: Formal Chaos}
and the phase transitions they induce in this more general setting.
That said, it is clear that such results do not consist of straightforward extensions of Theorem \ref{Theorem: Main}.

Firstly, to our knowledge, \eqref{Equation: Assumption} are the most general conditions
known to ensure that the moments of $\msf Q_\ka(t)$ are finite (at least for small enough times) under Assumption \ref{Assumption}; see \cite[Theorem 3.14-(iii)
and Lemma 4.4]{ChenDeyaOuyangTindel}.
Beyond that,
different proof techniques must be used to study the deterministic
and fluctuations terms in \eqref{Equation: Formal Chaos}.

Secondly, if $H_i\geq1/2$
is not satisfied, then $\msf Q_\ka(t)$'s moments may still be finite
(see \cite{ChenRoughSpace,ChenRoughTime,ChenDeyaOuyangTindel}),
but the moment formulas involved become much more difficult to analyze.
As highlighted by \eqref{Equation: g} and \eqref{Equation: h},
at a purely technical level, the issues are
that the covariance kernels of rough fractional noises are not positive measures
and that their Fourier transforms do not vanish at infinity.
However, these obstacles are not merely technical:
The increased roughness of the noise can induce entirely new interactions between
$\partial D$'s geometry and $\xi$'s statistical properties, which require a dedicated
analysis.

In order to illustrate this, consider
the constant $\mf c$ in \eqref{Equation: c Constant}. When $H_j<1/2$, the covariance
kernel $\msf g$ cannot be written as a product of Riesz kernels and Dirac masses.
Instead, one must use the spectral measure of $\xi$, which corresponds to $(\msf h\otimes\msf g)$'s
Fourier transform. In particular, the constant $\mf c$ in that setting is written as
\begin{align}
\label{Equation: Rough c Constant}
\mf c=\left(\prod_{j=1}^d\frac{\Ga(2H_j+1)}{2\Ga(H_j)\Ga(1-H_j)}\int_{\mbb R^d}|\hat{\bs 1_{D}}(\om)|^2\prod_{j=1}^d|\om_j|^{1-2H_j}\d\om\right)^{1/2},
\end{align}
where $\hat{\mbf 1_D}$ denotes the Fourier transform of $D$'s indicator function.
The problem with this expression is that if $H_j<1/2$ for some $j\geq1$,
then the boundedness of $D$ is no longer sufficient to guarantee the finiteness of $\mf c$:
Given that $|\om_j|^{1-2H_j}$ diverges as $\om_j\to\pm\infty$ when $H_j<1/2$,
the finiteness of $\mf c$ requires a quantitative control on the decay rate of $|\hat{\bs 1_{D}}(\om)|^2$,
the latter of which depends on the regularity of $\partial D$.
When that condition fails, one expects a different constant and scaling of $t$ in \eqref{Equation: Variance}.

In summary, given that these cases fall outside the scope of the methods used in this paper,
we leave their analysis open for future work.
In particular, we refer to Section \ref{Section: Open Problems} for a more detailed
discussion of conjectures and open problems related to extending 
Theorem \ref{Theorem: Main} and the resulting phase transitions beyond Assumption \ref{Assumption}.

\subsection{Pioneering Works and the Anderson Hamiltonian}

The study of small-time asymptotics
of diffusion problems perturbed by singular noises was, to the best of
our knowledge, pioneered in \cite{Mouzard}. Therein, Mouzard constructed the
two-dimensional Anderson Hamiltonian (AH)
\[\mc H_\ka=-\tfrac12\De+\ka\xi\]
on compact two-dimensional smooth Riemannian manifolds, assuming $\xi$ is
a time-independent white noise (i.e., $H_0=1$ and $H_1=H_2=1/2$). Then, he proved a first-order Weyl law for $\mc H_\ka$'s
eigenvalue counting function. In view of standard Abelian-Tauberian theorems, this can be recast
as the asymptotic $\mr{Tr}[\mr e^{-t\mc H_\ka}]\sim\mr{Tr}[\mr e^{-t\mc H_0}]$. See also
\cite{BailleulDangMouzard}, where this trace is studied directly for two-dimensional boundaryless manifolds.
Following up on this, Matsuda and van Zuijlen \cite{MatsudaVanZuijlen} extended Mouzard's Weyl law
to general time-independent singular noises in arbitrary dimensions,
as well as Neumann boundary conditions.

These results served as the main inspiration for
\cite{GaudreauLamarrePan}, as well as this paper.
We expect that our methods can also be applied to the AH with general time-independent fractional noise.
That being said,
since we have a particular interest in the
case where $\xi$ is time-dependent in this paper (due to the identification of $H_0$
as one of the three mechanisms controlling the total mass), we focus
our attention exclusively on the PAM  in the interest of
brevity.

\subsection{Organization and Outline of Main Technical Contributions}

The remainder of this paper is organized as follows:
In Section \ref{Section: Applications},
we explain our construction of $u_\ka$'s total mass for fixed $t$,
and we provide further context on the interpretation and applications of our main result.
Then, in Section \ref{Section: Open Problems}, we discuss conjectures and
open problems that naturally arise from our results.

In Section \ref{Section: Intersection Local Times}, we provide the main
technical tools used in our proofs, and take this opportunity to outline
some of the important ideas in our proof of Theorem \ref{Theorem: Main}
and our construction of the total mass (i.e., Proposition \ref{Proposition: L^2} below).
This section contains the first main technical input of the paper, namely,
Proposition \ref{Proposition: SILT Expectation Limits}. This result consists of
exact asymptotics for the expectation of the Brownian self-intersection-type
functionals associated with $\xi$'s covariance kernel. This is what allows us to
\begin{enumerate}
\item characterize when $\msf Q_\ka(t)$ needs a renormalization (i.e., $\rho\leq0$),
as well as the exact renormalization constants (i.e., Definition \ref{Definition: Renormalization} and Proposition \ref{Proposition: L^2}),
\item understand the precise relationship between the renormalization singularity and
the magnitude of $\msf Q_\ka(t)$'s expectation asymptotics when $\rho\leq0$ (i.e.,
Corollary \ref{Corollary: Renormalization} and Remark \ref{Remark: eps/t 1}), and
\item identify the exact constants in $\msf Q_\ka(t)$'s expectation asymptotics
(i.e., the constant $\mf b$ and the factor $\frac1{(1+\rho)\rho}$).
\end{enumerate}

Finally, using the tools from Section \ref{Section: Intersection Local Times}, in Section \ref{Section: L^2}--\ref{Section: Var}, we
respectively provide our construction of the total mass,
prove the expectation asymptotics in
Theorem \ref{Theorem: Main}, and then prove the variance asymptotics in
Theorem \ref{Theorem: Main}. The main novel technical input here comes from
Section \ref{Section: Var}, which provides an exact asymptotic
for the standard deviation. This is what allows us to characterize
how the space and time Hurst indices each contribute to the average
size of the fluctuations term in \eqref{Equation: Formal Chaos}.

\subsection*{Acknowledgments}

The authors thank Cyril Labb\'e for insightful conversations on
the contents of Section \ref{Section: Open Problems}.

\section{Geometry, Fluctuations, and Renormalization}
\label{Section: Applications}

Our main purpose in this section is to provide additional context
on two of the three mechanisms that control the formal total mass
expansion in \eqref{Equation: Formal Chaos}.
Namely, we discuss the geometry mechanism in
Section \ref{Section: Geometry}, and the
renormalization mechanism in 
Section
\ref{Section: Renormalization}.
Then, in Section \ref{Section: Transitions}, we state two phase
transition results, and we provide a phase diagram that summarizes
the meaningful interactions between $W$'s parameters.

\subsection{Geometry}
\label{Section: Geometry}

In the introduction, we hinted at the fact that the rate of heat diffusion through $\partial D$ can be quantified by geometric
properties of the latter.
This can be made concrete with the following results:

\begin{notation}
Given two nonvanishing functions $f$ and $g$ defined in a neighborhood of zero, we write
"$f(t)\lesssim g(t)$" if there exists $\th,C>0$ such that
$f(t)\leq C g(t)$ for all $t\in(0,\th)$, and we write
"$f(t)\asymp g(t)$" if $f(t)\lesssim g(t)$ and $g(t)\lesssim f(t)$.
\end{notation}

\begin{theorem}[{\cite[(1.16)]{BurchardSchmuckenschlager}}]
\label{Theorem: Heat Content Lower Limit}
For every bounded domain $D$,
\begin{align}
\label{Equation: Heat Content Lower Limit}
|D|-\int_D u_0(t,x)\d x\gtrsim t^{1/2}.
\end{align}
In other words, the total mass cannot converge to $|D|$ at a faster rate than $t^{1/2}$.
\end{theorem}

\begin{theorem}[\cite{VanDenBergPerimeter,VanDenBergDavies}]
\label{Theorem: Regular Heat Content}
If $\partial D$ is sufficiently regular
(such as the $R$-smooth condition in \cite[Definition 6.1]{VanDenBergDavies}), then
the lower bound \eqref{Equation: Heat Content Lower Limit} is optimal; in fact,
\begin{align}
\label{Equation: Regular Heat Content}
|D|-\int_D u_0(t,x)\d x\sim\frac{\sqrt{2}|\partial D|}{\sqrt{\pi}}t^{1/2},
\end{align}
where $|\partial D|$ is the Hausdorff measure of $D$'s boundary.
\end{theorem}

\begin{theorem}[\cite{VanDenBerg}]
\label{Theorem: Fractal Heat Content}
If there exists some $\msf m\in(d-1,d)$ such that
\[|\{x\in D:\mr{dist}(x,\partial D)<r\}|\asymp r^{d-\msf m},\qquad\text{as }r\to0\]
(which can happen if $\partial D$ is fractal)
and $\partial D$ satisfies a technical uniform capacity condition
(\cite[(1.9)]{VanDenBerg}), then the lower bound
\eqref{Equation: Heat Content Lower Limit} is no longer optimal; in fact,
\begin{align}
\label{Equation: Fractal Heat Content}
|D|-\int_D u_0(t,x)\d x\asymp t^{(d-\msf m)/2}.
\end{align}
\end{theorem}

\begin{remark}
If we assume more boundary regularity than the $R$-smooth condition of \cite{VanDenBergDavies}, then
\eqref{Equation: Regular Heat Content}
can be improved; see, e.g., \cite{VanDenBergLeGall} and \cite[Section 2.3]{Gilkey}.
\end{remark}

\subsection{Renormalization}
\label{Section: Renormalization}

We first outline our construction of the total mass, and then we explain the connection between the renormalization
and the expectation asymptotics in Theorem \ref{Theorem: Main}.

\subsubsection{Construction of the Total Mass}

For every $\eps>0$,
define the mollifier
\begin{align}
\label{Equation: Mollifier}
f_\eps(t,x)=\frac{\eps/2}{\pi(t^2+(\eps/2)^2)}\frac{\mr e^{-|x|^2/\eps}}{(\pi\eps)^{d/2}},\qquad (t,x)\in\mbb R\times\mbb R^d.
\end{align}
We then define the mollified noise
\[\xi_\eps(t,x)=\xi*f_{\eps}(t,x)=\int_{\mbb R\times\mbb R^d}f_{\eps}(t-s,x-y)\d \xi(s,y),\]
where we use $*$ to denote the convolution.
Given that the function $f_\eps$ is even,
we note that $\xi_\eps$ is a centered Gaussian process with covariance kernel
\begin{align}
\label{Equation: hg eps}
(\msf h\otimes\msf g)_\eps=(\msf h\otimes\msf g)*f_\eps*f_\eps=(\msf h\otimes\msf g)*f_{2\eps}
\end{align}
($f_\eps*f_\eps=f_{2\eps}$ because the Cauchy and Gaussian distributions are both stable).
More generally, for any $s,t\in\mbb R$, $x,y\in\mbb R^d$, and $\eps_1,\eps_2>0$,
one has
\begin{align}
\label{Equation: Covariance for eps1 and eps2}
\mbf E\big[\xi_{\eps_1}(s,x)\xi_{\eps_2}(t,y)\big]=(\msf h\otimes\msf g)_{(\eps_1+\eps_2)/2}(s-t,x-y).
\end{align}
In particular, $\xi_\eps$ has continuous sample paths almost surely.
With this in hand, we now consider
$u_{\ka,\eps}$, which is the formal solution of the PDE
\begin{align}
\label{Equation: u epsilon}
\partial_tu_{\ka,\eps}(t,x)=\big(\tfrac12\De+\ka\xi_\eps(t,x)\big)u_{\ka,\eps}(t,x),\qquad t>0\text{ and }x\in D,
\end{align}
with the same boundary and initial conditions as $u_\ka$.
In rigorous terms, we define the total mass of \eqref{Equation: u epsilon} using a probabilistic/Feynman-Kac formulation:

\begin{definition}
\label{Definition: Approximate Total Mass}
Let $B$ be a standard Brownian motion on $\mbb R^d$ (which is started at
the origin). For every $x\in\mbb R^d$, we define
\begin{align}
\label{Equation: Brownian Coupling}
B^x=x+B.
\end{align}
For any stochastic process $Z$ (e.g., $B^x$), we denote the hitting time
\[\tau_D(Z)=\inf\{s\geq0:Z(s)\not\in D\}.\]
Then, we define the total mass of $u_{\ka,\eps}$ (which we denote $\msf M_{\ka,\eps}(t)$) as
\begin{align}
\label{Equation: Feynman-Kac}
\int_Du_{\ka,\eps}(t,x)\d x:=\msf M_{\ka,\eps}(t)=\int_D\mbf E_B\left[\mbf 1_{\{\tau_D(B^x)>t\}}\mr e^{\ka\int_0^t\xi_\eps(t-u,B^x(u))\d u}\right]\d x
\end{align}
for every $\eps,t>0$, where we assume that $\xi$ and $B$ are independent,
and $\mbf E_B$ means that we are only taking the expectation with respect to $B$---i.e., conditional on $\xi$.
(This is clearly well-defined since $\xi_\eps$ is continuous.)
\end{definition}

\begin{remark}
The justification for Definition \ref{Definition: Approximate Total Mass}
is that the infinitesimal generator of the semigroup
\[p^D_t\phi(x)=\mbf E_B\left[\mbf 1_{\{\tau_D(B^x)>t\}}\phi\big(B^x(t)\big)\right],\qquad t>0,~\phi\in L^2(D)\]
is the Dirichlet Laplacian, defined using the quadratic form $(\phi,\psi)\mapsto\frac12\langle\nabla\phi,\nabla\psi\rangle$
with form domain $H^1_0(D)$; see,
e.g., \cite[(3.14), (3.15), and Theorem 4.11]{Sznitman}.
In particular, when $\ka=0$, \eqref{Equation: Feynman-Kac} reduces
to the Feynman-Kac formula for the heat content 
\begin{align}
\label{Equation: Heat Content FK}
\int_Du_0(t,x)\d x=\int_D\mbf P\left[\tau_D(B^x)>t\right]\d x.
\end{align}
Then, the exponential weight in
\eqref{Equation: Feynman-Kac} is the standard way to introduce a pointwise multiplication
by $\xi_{\eps}(t,x)$ in the PDE's dynamics.
When $\partial D$ is sufficiently regular and the initial and boundary conditions
are equal on $\partial D$ (which is not the case in this paper), it is known that the Feynman-Kac formula provides the unique strong
solution to \eqref{Equation: u epsilon}; see, e.g., \cite[Section 6, Theorem 5.2]{Friedman}
for more details.
\end{remark}

The Stratonovich approach to construct the PAM is to define $u_\ka$ as the $\eps\to0$ limit of \eqref{Equation: u epsilon}.
However a renormalization is sometimes necessary for this limit to exist:

\begin{definition}
\label{Definition: Renormalization}
Define the renormalization function $\msf c:(0,\infty)\to\mbb R$ as follows:
\begin{enumerate}
\item If $\rho>0$, then no renormalization
is necessary; hence $\msf c(\eps)=0$ for all $\eps>0$.
\item If $\rho=0$, then $\xi$ needs a logarithmic
renormalization: If $\mf b$ is the constant defined in \eqref{Equation: b Constant}, then
\begin{align}
\label{Equation: Renormalization - Critical}
\msf c(\eps)=\mf b\,\log(1/\eps),\quad\eps>0.
\end{align}
\item If $\rho<0$, then $\xi$ requires a power-law renormalization:
\begin{align}
\label{Equation: Renormalization - Supercritical}
\msf c(\eps)=\mf d\,\eps^{\rho},\qquad\eps>0,
\end{align}
where we define the constant
\begin{align}
\label{Equation: Renormalization - Supercritical constant}
\mf d=
\frac{\Ga(-\rho)}{2^H\Ga(d-H)\Ga(1-H_0)}\prod_{j=0}^d\frac{\Ga(2H_j+1)}{\Ga(H_j)}
\int_0^{\pi/2}\frac{(\cos\theta)^{-(2H_0-1)}(\sin\theta)^{d-H}}{(\cos\theta+\sin\theta)^{-\rho}}\d\theta.
\end{align}
\end{enumerate}
\end{definition}

We may now construct the PAM's total mass and $\msf Q_\ka(t)$ as follows:

\begin{proposition}
\label{Proposition: L^2}
Under Assumption \ref{Assumption},
there exists a constant $\vartheta_\ka\in(0,\infty]$ such that for every $t\in(0,\vartheta_\ka)$,
there exists a square-integrable random variable $\msf M_\ka(t)$ satisfying
\[\lim_{\eps\to0}\mbf E\big[\big|\msf M_{\ka,\eps}(t)\mr e^{-t\ka^2\msf c(\eps)}-\msf M_\ka(t)\big|^2\big]=0.\]
\end{proposition}

\begin{definition}
\label{Definition: L^2}
For every $\ka>0$ and $t\in(0,\th_\ka)$, we define $u_\ka$'s total mass as
\[\int_Du_\ka(t,x)\d x:=\msf M_\ka(t).\]
We then define $\msf Q_\ka(t)$ as \eqref{Equation: Qkappa}.
\end{definition}

See Section \ref{Section: L^2} for a proof.

\begin{remark}
The proof of Proposition \ref{Proposition: L^2} relies on a standard
procedure, namely, a combination of the Feynman-Kac formula with
known estimates of Brownian intersection local times;
see \cite{ChenRoughSpace,ChenRoughTime,ChenDeyaOuyangTindel,HuHuangNualartTindel}. We nevertheless provide a complete proof since we are not aware of
a prior construction of the PAM (or its total mass) on general bounded domains. In any case, all the necessary ingredients
in the proof are also needed elsewhere.
We also note that
Stratonovich constructions are typically much more involved than what we do in Proposition
\ref{Proposition: L^2}, and use sophisticated analytic tools like paracontrolled calculus or regularity structures \cite{GIP,Hairer};
see, e.g., \cite{ChenDeyaOuyangTindel2,HairerKPZ,HairerLabbe}.
We are able to sidestep these difficulties entirely due to the availability of a probabilistic
representation for $u_{\ka,\eps}$, and the fact that we are only looking at a real-valued
observable of the solution (i.e., the total mass) at a fixed time, rather than the whole
random field solution
over $(t,x)\in[0,\infty)\times D$.
\end{remark}

\begin{remark}
The only situation where it is necessary to take a finite $\th_\ka$
in Proposition \ref{Proposition: L^2} is when $d-H=1$. Indeed,
these models are known to be "moment-critical" (see \cite[(1.8)]{ChenDeyaOuyangTindel}), in the sense that the PAM's
moments only exist up to some finite critical time. Given that we are only
interested in small-time asymptotics, however, this is not a distinction
that is important in this paper.
\end{remark}

\subsubsection{Renormalizations Determine Expectation Asymptotics}

With Definition \ref{Definition: Renormalization} in hand, we are now able to clarify the relationship between the
renormalization criticality and the expectation asymptotics in Theorem \ref{Theorem: Main}:

\begin{corollary}
\label{Corollary: Renormalization}
Suppose that Assumption \ref{Assumption} holds, and let $\ka>0$.
\begin{enumerate}
\item If $\rho>0$, then
$\displaystyle\big|\mbf E\big[\msf Q_\ka(t)\big]\big|=t\cdot o(1)$ (consistent with $\msf c(\eps)=0$).
\vspace{5pt}
\item If $\rho\leq0$, then
$\displaystyle\big|\mbf E\big[\msf Q_\ka(t)\big]\big|\asymp t\cdot\msf c(t)$ (except when $H_0=1/2$, due to $\mf b=0$).
\end{enumerate}
\end{corollary}

\begin{remark}
\label{Remark: eps/t 1}
Corollary \ref{Corollary: Renormalization} can largely be explained by
Brownian scaling and the parabolic scaling in the mollifier
\begin{align}
\label{Equation: Parabolic Scaling}
f_{\eps}(t,x)=\eps^{-d/2-1}f_1(t/\eps,x/\sqrt\eps).
\end{align}
As a result of this combination, in the proof of Proposition \ref{Proposition: L^2}, the singularity in $\msf M_{\ka,\eps}(t)$ involves a term that depends on the
ratio $\eps/t$. Consequently, any diverging $\eps$-dependent term that requires a renormalization
is accompanied by a "mirror image" $t$-dependent term, which survives the $\eps\to0$ limit.
See Remark \ref{Remark: eps/t 2} for more details.
\end{remark}

\subsection{Phase Transitions}
\label{Section: Transitions}

\subsubsection{The Geometry-Fluctuations Transition}

Under Assumption \ref{Assumption}, it is always the
case that $\rho>-1/2$ since $d-H<2H_0-1/2$.
In particular, $\mbf E\big[\msf Q_\ka(t)\big]=o(t^{1/2})$ as $t\to0$
by Theorem \ref{Theorem: Main}. Therefore, in the setting of this
paper, the geometry and deterministic contributions to \eqref{Equation: Formal Chaos}
do not intersect.

However, we see from \eqref{Equation: Variance} that
the fluctuations term can catch up to the heat content when $H_0=1/2$.
More specifically, we get the following dichotomy as an immediate consequence of Theorems \ref{Theorem: Main},
\ref{Theorem: Regular Heat Content}, and \ref{Theorem: Fractal Heat Content},
and a straightforward application of Chebyshev's inequality:

\begin{corollary}
\label{Corollary: Geometry-Fluctuations}
Suppose that Assumption \ref{Assumption} holds,
and let $\ka>0$.
\begin{enumerate}
\item If $H_0>1/2$ or $\partial D$ is sufficiently irregular, then the deterministic
heat diffusion mechanism dominates $u_\ka$'s total mass asymptotics. More specifically:
\begin{enumerate}
\item If $H_0>1/2$ and \eqref{Equation: Regular Heat Content} holds, then
\[\lim_{t\to0} t^{-1/2}\big(|D|-\msf M_\ka(t)\big)=\frac{\sqrt{2}|\partial D|}{\sqrt{\pi}}\qquad\text{in probability}.\]
\item If \eqref{Equation: Fractal Heat Content} holds for any $\msf m\in(d-1,d)$, then
\[\lim_{t\to0}\frac{1}{\log t}\log\big||D|-\msf M_\ka(t)\big|=\frac{d-\msf m}{2}\qquad\text{in probability}.\]
\end{enumerate}
\item If $H_0=1/2$ and \eqref{Equation: Regular Heat Content} holds, then as $t\to0$, one has
\[\mbf E\big[|D|-\msf M_\ka(t)\big]\sim\frac{\sqrt{2}|\partial D|}{\sqrt{\pi}}t^{1/2}
\quad\text{and}\quad
\sqrt{\mbf{Var}\big[|D|-\msf M_\ka(t)\big]}\sim\ka\mf c\,t^{1/2}.\]
In particular, the noise now has an impact on the first-order total mass asymptotic via $\ka$ and the constant $\mf c$ in \eqref{Equation: c Constant}.
\end{enumerate}
\end{corollary}

\subsubsection{The Fluctuations-Renormalization Transition}

When $\rho\leq0$, there is a possibility for the expectation asymptotic
in Theorem \ref{Theorem: Main} to either be the same size or larger
than the standard deviation. This leads to the following phase transition,
which is an immediate consequence of
Theorem \ref{Theorem: Main},
the definition of $\mf b$ in \eqref{Equation: b Constant},
and a straightforward application of Chebyshev's inequality.

\begin{corollary}
\label{Corollary: Fluctuations-Renormalization}
Suppose that Assumption \ref{Assumption} holds,
and let $\ka>0$.
\begin{enumerate}
\item If $\rho>H_0-1$, then
\[\mbf E\big[\msf Q_\ka(t)\big]=o(t^{H_0})
\qquad\text{and}\qquad
\mbf E[\msf Q_\ka(t)^2]^{1/2}\sim\ka\mf c\,t^{H_0}.\]
\item If $\rho=H_0-1$ and $H_0\neq1$, then
\[\mbf E\big[\msf Q_\ka(t)\big]\sim\frac{\ka^2|D|(2H_0-1)}{(H_0-1)2^{H}}\prod_{j=1}^d\frac{\Ga(2H_j+1)}{\Ga(H_j)}\,t^{H_0}
\qquad\text{and}\qquad
\sqrt{\mbf{Var}[\msf Q_\ka(t)]}\sim\ka\mf c\,t^{H_0}.\]
\item Otherwise:
\begin{enumerate}
\item If $\rho=H_0-1$ and $H_0=1$, then
\[\lim_{t\to0}\frac{\msf Q_\ka(t)}{t\log(1/t)}=-\frac{\ka^2|D|}{2^{d-1}}\prod_{j=1}^d\frac{\Ga(2H_j+1)}{\Ga(H_j)}\qquad\text{in probability.}\] 
\item If $\rho<H_0-1$, then
\[\lim_{t\to0}\frac{\msf Q_\ka(t)}{t^{1+\rho}}
=\frac{\ka^2|D|H_0(2H_0-1)}{2^H(1+\rho)\rho}\prod_{j=1}^d\frac{\Ga(2H_j+1)}{\Ga(H_j)}\qquad\text{in probability.}\]
\end{enumerate}
\end{enumerate}
\end{corollary}

\begin{remark}
There is a vast literature on statistical inference for parameter-dependent
stochastic PDEs; see, e.g., the surveys \cite{Stats4SPDEs1,Stats4SPDEs2}, as well as the comprehensive
bibliography maintained on the
website \url{https://sites.google.com/view/stats4spdes/bibliography}.
To the best of our knowledge, none of these results concern renormalized
Stratonovich equations. In this context, it is interesting to note that Corollary \ref{Corollary: Fluctuations-Renormalization}-(3)
suggests that it is possible to recover (in the sense of a consistent estimator) some information
about $\ka$ and the parameters of $W$ from the total mass' small-$t$ asymptotics.
We leave the question of whether this can be leveraged into meaningful statistical inference procedures open.
\end{remark}

\subsubsection{Summary and Diagram}

In Figure \ref{Figure: Phase Diagram} below,
we provide a complete phase diagram for the transitions that emerge from Theorem \ref{Theorem: Main}.
\begin{figure}[htbp]
\begin{center}
\includegraphics[width=\textwidth]{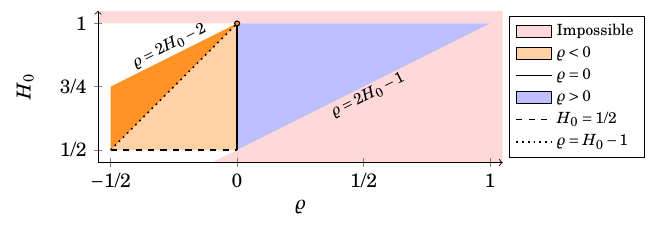}
\caption{Complete Phase Diagram for Assumption \ref{Assumption}.}
\label{Figure: Phase Diagram}
\end{center}
\end{figure}
The region shaded in red represents impossible combinations of $\rho$
and $H_0$ (i.e., $H_0>1$ or $\rho\geq2H_0-1$, the latter of which implies $H\geq d$).
The other shaded regions are permissible under Assumption \ref{Assumption}
(however, when $d=1$ there are further restrictions). In particular, the diagonal edge on the line $\rho=2H_0-2$ comes
from the assumption that $d-H\leq1$, and the lower bound $\rho>-1/2$ comes from the assumption that $d-H<2H_0-1/2$.

The blue region represents the regime $\rho>0$ wherein no renormalization is required.
The line $\rho=0$ and the orange region $\rho<0$ are respectively the logarithmic and power-law renormalized regimes,
as per Definition \ref{Definition: Renormalization}.
This visualization makes
it clear that the phase transition phenomena uncovered in this paper are exclusive
to the renormalized models.

More specifically: The dashed line at $H_0=1/2$ represents both the vanishing expectation criticality
in Theorem \ref{Theorem: Main}-(3), and the critical threshold of the
geometry-fluctuations transition in Corollary \ref{Corollary: Geometry-Fluctuations}.
The dotted line at $\rho=H_0-1$ represents the critical threshold of the
fluctuations-renormalization transition in Corollary \ref{Corollary: Fluctuations-Renormalization}.
The dark orange dot and the
dark orange region
represent the cases where $\msf Q_\ka(t)$'s expectation decays more slowly
than its standard deviation, which gives rise to the consistent estimators stated in Corollary \ref{Corollary: Fluctuations-Renormalization}-(3).

\section{Conjectures and Open Problems}
\label{Section: Open Problems}

In principle, the parameter $H_0$ could take any
value in $(0,1]$, and the parameter $\rho$ could take
any value in $(-d-1,1)$ (though not every combination of
$H_0$ and $\rho$ in that rectangle is possible,
as per the red regions in Figure \ref{Figure: Phase Diagram}). It is natural to wonder what becomes of Theorem
\ref{Theorem: Main}, the three competing mechanisms
in \eqref{Equation: Formal Chaos}, and the resulting
phase transitions illustrated in Figure \ref{Figure: Phase Diagram} beyond the cases that satisfy Assumption \ref{Assumption}. In this section, we discuss various conjectures and open problems related to these extensions.

\subsection{Rough Time and Fluctuations}

Arguably the simplest extension of our setting
would be in the case where $H_0<1/2$, $H_j\geq1/2$ for all $j\geq1$,
and the moments of $\msf Q_\ka(t)$ are still finite.

\begin{conjecture}
\label{Conjecture: Finite Moments with Rough Time}
If $H_0<1/2$, $H_j\geq1/2$ for all $j\geq1$,
and the moments of $\msf Q_\ka(t)$ are finite for small $t>0$, then
the statement of Theorem \ref{Theorem: Main} holds without
modification.
\end{conjecture}

If true, then Conjecture \ref{Conjecture: Finite Moments with Rough Time} would extend
the geometry-fluctuations phase transition in Corollary \ref{Corollary: Geometry-Fluctuations},
into a supercritical regime where the fluctuations are the sole dominant term in the total mass asymptotics.
More specifically, conjecture \ref{Conjecture: Finite Moments with Rough Time} suggests that if either
\begin{enumerate}
\item $H_0<1/2$ and \eqref{Equation: Regular Heat Content} hold, or
\item $H_0<(d-\msf m)/2$ and \eqref{Equation: Fractal Heat Content} hold
for some $\msf m\in(d-1,d)$,
\end{enumerate}
then, as $t\to0$, one has
\[\mbf E\big[|D|-\msf M_\ka(t)\big]=o(t^{H_0})
\quad\text{and}\quad
\sqrt{\mbf{Var}\big[|D|-\msf M_\ka(t)\big]}\sim\ka\mf c\,t^{H_0}.\]

Corollary \ref{Corollary: Geometry-Fluctuations}
and Conjecture \ref{Conjecture: Finite Moments with Rough Time}
also motivate the following problem:

\begin{openproblem}
\label{Open Problem: Fluctuations}
For every $n\geq1$, define the stochastic process
\[X_n(t)=n^{H_0}\big(\msf Q_\ka(t/n)-\mbf E[\msf Q_\ka(t/n)]\big),\qquad 0\leq t<n\th_\ka.\]
Does there exist a stochastic process $X=(X_t:t\geq0)$ such that $X_n\to X$ in distribution
as $n\to\infty$?
If so, then what is $X$'s distribution?
\end{openproblem}

Indeed, understanding $X$ would help explain the mechanism responsible
for $u_\ka$'s total mass asymptotics when the geometry-fluctuations threshold
of $H_0=1/2$ is either attained or crossed.

\subsection{Finite Moments with Rough Space}

Next, suppose that in addition to $H_0$, we also allow some of the space Hurst indices to be
smaller than $1/2$. In this situation, the known conditions that ensure that
$\msf Q_\ka(t)$'s moments are finite are more complicated than
\eqref{Equation: Assumption}; see \cite[(1.13)]{ChenRoughSpace},
\cite[(1.6)]{ChenRoughTime}, and
\cite[(1.7) and (1.8)]{ChenDeyaOuyangTindel}.
Nevertheless, in that case we expect the following:

\begin{conjecture}
\label{Conjecture: Finite Moments with Rough Space}
Suppose that $H_j<1/2$ for at least one $j\geq1$,
and that $\msf Q_\ka(t)$'s moments are finite for small enough $t>0$.
If the constant $\mf c$ in \eqref{Equation: Rough c Constant} is finite,
then the statement of Theorem \ref{Theorem: Main} holds without
modification.
\end{conjecture}

We expect that the real difficulty in this setting will be as follows:

\begin{openproblem}
Characterize when $\mf c<\infty$. Then, understand how the constant and the power of
$t$ in \eqref{Equation: Variance} must be modified when $\mf c=\infty$.
\end{openproblem}

Thanks to \eqref{Equation: Rough c Constant},
we expect that the constant and the power of $t$ that appear in the variance asymptotic
when $\mf c=\infty$ depend on some combination of
$H_0$, the growth rates of the functions $\om_j\mapsto|\om_j|^{1-2H_j}$ for all $j$'s
such that $H_j<1/2$, and the decay rate of the Fourier transform $|\hat{\bs 1_{D}}(\om)|^2$
as $|\om_j|\to\infty$. This suggests that, in sharp contrast to our setting,
the fluctuations term is also influenced by $\partial D$'s geometry when
the spatial covariance is rough---not just $H_0$. If true, then this could significantly alter the critical thresholds for
both the geometry-fluctuations and fluctuations-renormalization transitions in Corollaries \ref{Corollary: Geometry-Fluctuations}
and \ref{Corollary: Fluctuations-Renormalization}.

\subsection{Beyond Finite Moments}

Lastly, we consider the case when $\msf Q_\ka(t)$'s moments blow up for all $t>0$.
Without the expectation and standard deviation, one must find alternate
means of estimating the deterministic and fluctuations terms in \eqref{Equation: Formal Chaos}.
Nevertheless, we expect that the phase transitions uncovered in this paper should still occur, and that entirely new transitions should emerge.

\subsubsection{Renormalization Thresholds and Deterministic Asymptotics}

It is known that increasing the singularity of the noise
can lead to renormalization criticalities beyond those stated
in Definition \ref{Definition: Renormalization}. For instance,
if $\xi$ is a time-independent white noise on $\mbb R^3$, then
it was shown in \cite{HairerLabbe} that the PAM needs two renormalizations:
A power law of order $\eps^{-1/2}$
and a smaller correction of order $\log(1/\eps)$.

\begin{remark}
In \cite[(1.1)]{HairerLabbe}, the power law renormalization for the three-dimensional PAM
is stated as $\eps^{-1}$ instead of $\eps^{-1/2}$. However, this is only because the mollifier is
scaled differently in that paper; see \cite[paragraph preceding (E$_\eps$)]{HairerLabbe}.
\end{remark}

We expect that these additional criticalities should induce new phase transitions
in the asymptotics of the deterministic part of \eqref{Equation: Formal Chaos}.
Using \cite[(1.1)]{HairerLabbe} as a guide, we expect the following:

\begin{conjecture}
\label{Conjecture: More Renormalization Transitions}
If $d=3$, $H_0=1$, and $H_1=H_2=H_3=1/2$, then
there exist constants $\mf A,\mf B\neq0$ such that
\[\msf Q_\ka(t)=\ka^2\mf A\,t^{1/2}+\ka^4\mf B\,t\log(1/t)+O_{\mbf P}(t)\qquad\text{as }t\to0.\]
More generally, if a model of the PAM requires multiple renormalization functions of decreasing order, say $\msf c_1(\eps),\ldots,\msf c_n(\eps)$, then there are corresponding
deterministic contributions of order $t\msf c_k(t)$ in the small-$t$ asymptotics of $\msf Q_\ka(t)$.
\end{conjecture}

\begin{remark}
If true, then Conjecture \ref{Conjecture: More Renormalization Transitions} would
also deepen the phase transition in Corollary \ref{Corollary: Fluctuations-Renormalization}
(i.e., fluctuations-renormalization), in the sense that there could
be multiple deterministic contributions of greater order than the fluctuations term.
\end{remark}

\subsubsection{A Geometric-Renormalization Transition}

As $\rho$ approaches $-1/2$ from above,
the power law $t^{1+\rho}$ in $\msf Q_\ka(t)$'s expectation asymptotic approaches $t^{1/2}$.
This gives rise to the tantalizing possibility
of a third phase transition, wherein $u_\ka$'s total mass asymptotics
can be dominated by a deterministic mechanism that comes from the noise;
presumably the terms $\msf D_{2n}(t)$ in \eqref{Equation: Formal Chaos}.
Thus:

\begin{openproblem}
If $\rho\leq-1/2$ (i.e., $d-H\geq2H_0-1/2$), does there sometimes exist a scaling function $\msf s$ such that
$\msf s(t)\big(|D|-\msf M_\ka(t)\big)$
converges in probability to a $\ka$-dependent constant as $t\to0$?
\end{openproblem}

A naive interpolation of the expectation asymptotic when $\rho>-1/2$
suggests that we could take $\msf s(t)=t^{-1-\rho}$.
One natural test case for this heuristic would be the three-dimensional time-independent white noise PAM,
which satisfies $\rho=-1/2$. In that case,
the prediction $\mathsf s(t)=t^{-1/2}$ is consistent with Conjecture \ref{Conjecture: More Renormalization Transitions}.

\section{Intersection Local Times and Main Proof Architecture}
\label{Section: Intersection Local Times}

In this section, we provide the main technical inputs of our proof technique,
and we outline how each
of these elements contributes to the proof of our main results.

\subsection{Main Statements}
\label{Section: Intersection Local Times Main}

\subsubsection{Feynman-Kac Formulas for $\eps>0$}

Our first tool consists of Feynman-Kac formulas for
$\msf M_{\ka,\eps}$'s moments when $\eps>0$. We begin with some notations
and definitions:

\begin{notation}
We use $B_1$ and $B_2$ to denote
i.i.d. copies of $B$, and adopt the same coupling
as in \eqref{Equation: Brownian Coupling} for $B^{x_i}_i$
(where $x_i\in\mbb R^d$).
\end{notation}

\begin{definition}
\label{Definition: Approximate SILT and MILT}
Let $\eps,t>0$ be fixed. Recall the definition of
$(\msf h\otimes\msf g)_{\eps}$ in \eqref{Equation: hg eps}.
We define $B$'s
approximate self-intersection local time (SILT) on $[0,t]$ as
\[\be^\eps_t(B)=
\int_{[0,t]^2}(\msf h\otimes\msf g)_\eps\big(s_2-s_1,B(s_2)-B(s_1)\big)\d s.\]
For any $x=(x_1,x_2)\in(\mbb R^d)^2$, we define $B^{x_i}_i$'s approximate mutual intersection local time
(MILT) on $[0,t]$ as
\[\al^\eps_t(x)=
\int_{[0,t]^2}(\msf h\otimes\msf g)_\eps\big(s_2-s_1,B^{x_2}_2(s_2)-B^{x_1}_1(s_1)\big)\d s.\]
\end{definition}

We are now in a position to state the Feynman-Kac formulas
that form the basis of our proof technique:

\begin{proposition}
\label{Proposition: Finite Epsilon Feynman-Kac}
Let $\ka\geq0$ and $t,\eps>0$ be fixed.
One has
\begin{align}
\label{Equation: Finite Epsilon Feynman-Kac E}
\mbf E\big[\msf M_{\ka,\eps}(t)\big]
&=\int_D\mbf E\left[\mbf 1_{\{\tau_D(B^x)>t\}}\mr e^{\frac{\ka^2}{2}\be^\eps_t(B)}\right]\d x,\\
\label{Equation: Finite Epsilon Feynman-Kac Var}
\mbf{Var}\big[\msf M_{\ka,\eps}(t)\big]
&=\int_{D^2}\mbf E\Bigg[\mbf 1_{\cap_{i\leq2}\{\tau_D(B^{x_i}_i)>t\}}\mr e^{\sum_{i=1}^2\frac{\ka^2}{2}\be^{\eps}_t(B_i)}
\left(\mr e^{\ka^2\al^\eps_t(x)}-1\right)\Bigg]\d x.
\end{align}
\end{proposition}

Proposition \ref{Proposition: Finite Epsilon Feynman-Kac} is a straightforward consequence of the
classical Feynman-Kac formula; see Section
\ref{Section: Proof of Proposition: Finite Epsilon Feynman-Kac} for a proof.

\subsubsection{Limits and Uniform Integrability}

In order to construct $\msf M_\ka(t)$
as an $L^2$ limit (i.e., Proposition \ref{Proposition: L^2}) and
derive Feynman-Kac formulas for $\msf M_\ka(t)$'s moments,
we need to understand the $\eps\to0$ limits of the approximate SILT's
and MILT's appearing in Proposition \ref{Proposition: Finite Epsilon Feynman-Kac},
as well as the uniform integrability properties of the prelimit variables.
For this purpose, we have the following results:

\begin{proposition}
\label{Proposition: Limit and UI of SILT}
Suppose that Assumption \ref{Assumption} holds.
\begin{enumerate}
\item If $\rho>0$, for every $t>0$, there exists a random variable $\ga_t(B)$ such that
\begin{align}
\label{Equation: subcritical SILT limit}
\lim_{\eps\to0}\be^\eps_t(B)=\ga_t(B)
\qquad\text{in $L^2$}.
\end{align}
Moreover, for every $c,t>0$, one has
\begin{align}
\label{Equation: subcritical SILT UI}
\sup_{\eps>0}\mbf E\left[\exp\big(c\be^\eps_t(B)\big)\right]\leq\mbf E\left[\exp\big(c\ga_t(B)\big)\right]<\infty.
\end{align}
\item If $\rho\leq0$, for every $t>0$, there exists a random variable $\ga_t(B)$
such that
\begin{align}
\label{Equation: renormalized SILT limit}
\lim_{\eps\to0}\Big(\be^\eps_t(B)-\mbf E\big[\be^\eps_t(B)\big]\Big)=\ga_t(B)
\qquad\text{in $L^2$}.
\end{align}
Moreover, for every $c>0$, there exists $\th_c>0$ such that for every $t\in(0,\th_c)$,
\begin{align}
\label{Equation: renormalized SILT UI}
\sup_{\eps>0}\mbf E\left[\exp\Big(c\big(\be^\eps_t(B)-\mbf E\big[\be^\eps_t(B)\big]\big)\Big)\right]<\infty
\quad\text{and}\quad
\mbf E\left[\exp\big(c\ga_t(B)\big)\right]<\infty.
\end{align}
\end{enumerate}
\end{proposition}

\begin{remark}
When $\rho>0$, it is customary to denote the limit of $\be^\eps_t(B)$
as $\be_t(B)$. $\ga_t(B)$ is instead more common as the notation for
the renormalized SILT, wherein we subtract $\be^\eps_t(B)$'s expectation. We adopt this nonstandard
notation in order to streamline our presentation.
\end{remark}

\begin{proposition}
\label{Proposition: Limit and UI of MILT}
Suppose that Assumption \ref{Assumption} holds.
For every $t>0$ and $x\in(\mbb R^d)^2$,
there exists a random variable $\al_t(x)$ such that
\begin{align}
\label{Equation: MILT Limit}
\lim_{\eps\to0}\al^\eps_t(x)=\al_t(x)
\qquad\text{in $L^p$ for all $p\geq1$}.
\end{align}
Moreover, for every $c>0$, there exists $\th_c>0$
such that for every $t\in(0,\th_c)$,
\begin{align}
\label{Equation: MILT UI}
\sup_{\eps>0,~x\in(\mbb R^d)^2}\mbf E\Big[\exp\big(c\al^\eps_t(x)\big)\Big],\sup_{x\in(\mbb R^d)^2}\mbf E\Big[\exp\big(c\al_t(x)\big)\Big]<\infty.
\end{align}
\end{proposition}

Propositions \ref{Proposition: Limit and UI of SILT} and \ref{Proposition: Limit and UI of MILT}
are essentially well-known---the tools we need consist of
standard material in the monograph \cite{ChenBook} and estimates found in
\cite{ChenRoughSpace,ChenDeyaOuyangTindel}. However, given that we did not find these exact
statements in the literature for all cases of Assumption \ref{Assumption}, we nevertheless
provide a proof in Section \ref{Section: Proof of Limit and UI of SILT and MILT}.

\subsubsection{Renormalization via Expectation Limits}

With Propositions \ref{Proposition: Limit and UI of SILT} and \ref{Proposition: Limit and UI of MILT}
in hand, it is clear that the limits of the Feynman-Kac formulas in
\eqref{Equation: Finite Epsilon Feynman-Kac E} and \eqref{Equation: Finite Epsilon Feynman-Kac Var}
will involve $\ga_t(B_i)$ and $\al_t(x)$. However, there is still one more crucial element that we need
to calculate these limits, which is coming from the renormalization in \eqref{Equation: renormalized SILT limit},
i.e., the subtraction of $\mbf E\big[\be^\eps_t(B)\big]$. As it turns out, this is the mechanism at the
origin of the renormalization functions in Definition \ref{Definition: Renormalization}, as is made clear
by the following scaling property and asymptotics:

\begin{lemma}
\label{Lemma: Scaling of Beta and Alpha}
For every $t,\eps>0$,
\[\be^\eps_t(B)\deq t^{1+\rho} \be^{\eps/t}_1(B).\]
Moreover, for every $x\in(\mbb R^d)^2$,
\[\al^\eps_t(x)
\deq
t^{1+\rho}
\al^{\eps/t}_1(x/\sqrt t).\]
In particular, taking $\eps\to0$, we get that
\begin{align}
\label{Equation: Gamma Alpha Scalings}
\ga_t(B)\deq t^{1+\rho}\ga_1(B)
\qquad\text{and}\qquad
\al_t(x)
\deq
t^{1+\rho}
\al_1(x/\sqrt t).
\end{align}
\end{lemma}

Lemma \ref{Lemma: Scaling of Beta and Alpha} follows from a straightforward combination
of the parabolic scaling \eqref{Equation: Parabolic Scaling} and Brownian scaling.
See Section \ref{Section: Proof of Scaling of Beta and Alpha} for a proof.

\begin{proposition}
\label{Proposition: SILT Expectation Limits}
Suppose that Assumption \ref{Assumption} holds, and
let $\ka>0$ be fixed.
\begin{enumerate}
\item If $\rho>0$, then
\begin{align}
\label{Equation: Separate Expectation Limits - Subcritical}
\lim_{\th\to0}\mbf E\big[\tfrac{\ka^2}{2}\be^\th_1(B)\big]=\frac{\ka^2\mf b}{(1+\rho)\rho}.
\end{align}
\item If $\rho=0$, then
there exist bounded and continuous functions $\msf b_i:[0,\infty)\to\mbb R$, $(i=1,2)$,
such that for every $\th>0$, one has
\begin{align}
\label{Equation: Separate Expectation Limits - Critical}
\mbf E\big[\tfrac{\ka^2}{2}\be^\th_1(B)\big]
=
\begin{cases}
\ka^2\mf b\,\ell(\th)+\ka^2\msf b_1(\th)&H_0<1,\\
\ka^2\mf b\,(1+\th)\log(1+1/\th)+\ka^2\msf b_2(\th)&H_0=1,
\end{cases}
\end{align}
where
\begin{align}
\label{Equation: ell function}
\ell(\th)=\frac{1}{2\Ga(H_0)\Ga(1-H_0)}\int_0^{2\pi}\left|\frac{\sin\theta}{\cos\theta}\right|^{2H_0-1}
\Ga_0\Big(\th\big(|\sin\theta|+|\cos\theta|\big)\Big)\d\theta,
\end{align}
and $\Gamma_0$ denotes the zero-parameter upper incomplete Gamma function, i.e.,
\[\Gamma_0(\th)=\int_{\th}^\infty \mr e^{-z}z^{-1}\d z,\qquad \th>0.\]
\item 
If $\rho<0$,
then there exists a bounded and continuous function $\msf z:(0,\infty)\to\mbb R$ such that
for every $\th>0$,
\begin{align}
\label{Equation: Separate Expectation Limits - Supercritical}
\mbf E\big[\tfrac{\ka^2}{2}\be^\th_1(B)\big]=\ka^2\mf d\,\th^{\rho}+\,\ka^2\msf z(\th).
\end{align}
Moreover,
$\msf z$ satisfies the limit
\begin{align}
\label{Equation: Separate Expectation Limits - Supercritical z limit}
\lim_{\th\to0}\msf z(\th)=
\begin{cases}
\frac{\mf b}{(1+\rho)\rho}&H_0>1/2,\\
0&H_0=1/2.
\end{cases}
\end{align}
\end{enumerate}
\end{proposition}

See Section \ref{Section: SILT Expectation Limits} for a proof. We may now anticipate
both the renormalization function in Definition \ref{Definition: Renormalization} and the
appearance of $\frac{\ka^2\mf b}{(1+\rho)\rho}$
in the expectation asymptotics using the following:

\begin{corollary}
\label{Corollary: Limits of E minus Renormalization}
Let $\ka>0$ be fixed.
\begin{enumerate}
\item If $\rho=0$, then there exists a constant $\mf e_\ka\in\mbb R$ such that
\[\lim_{\eps\to0}\Big(\mbf E\big[\tfrac{\ka^2}{2}\be^\eps_t(B)\big]-t\ka^2\msf c(\eps)\Big)=-\ka^2\mf b\,t\log(1/t)+\mf e_\ka\,t.\]
\item If $\rho<0$, then
\[\lim_{\eps\to0}\Big(\mbf E\big[\tfrac{\ka^2}{2}\be^\eps_t(B)\big]-t\ka^2\msf c(\eps)\Big)=\tfrac{\ka^2\mf b}{(1+\rho)\rho}\, t^{1+\rho}.\]
\end{enumerate}
\end{corollary}

See Section \ref{Section: Limits of E minus Renormalization}
for a proof.

\begin{remark}
\label{Remark: eps/t 2}
Following up on Remark \ref{Remark: eps/t 1},
we now have additional context on the mechanism driving Corollary \ref{Corollary: Renormalization}. More specifically, when $\rho=0$, the $t\log(1/t)$ term emerges from a term of the form $t\log(t/\eps)$,
and only the $t\log(1/\eps)$ requires a renormalization. When $\rho<0$, the emergence of
the $t^{1+\rho}$ is a bit more subtle, and relies on the combination of the scaling in Lemma \ref{Lemma: Scaling of Beta and Alpha} with the subsequent cancellations in \eqref{Equation: Separate Expectation Limits - Supercritical}
and the limit of $\msf z$ in \eqref{Equation: Separate Expectation Limits - Supercritical z limit}.
\end{remark}

We finish the statements of our technical results with
the following, which is a consequence of
the $H_0=1/2$ case of \eqref{Equation: Separate Expectation Limits - Supercritical} and \eqref{Equation: Separate Expectation Limits - Supercritical z limit}:

\begin{corollary}
\label{Corollary: SILT=0 when H_0=1/2}
Let $\rho<0$, and let $\ga_t(B)$ be the
random variable constructed in \eqref{Equation: renormalized SILT limit}. If $H_0=1/2$,
then $\ga_t(B)=0$ almost surely.
\end{corollary}

See Section \ref{Section: SILT=0 when H_0=1/2}
for the proof.

\subsection{Fourier Analysis}

Before embarking on the proofs of the results in this section, 
we introduce some Fourier analysis, which makes calculations easier.
To begin with, we use the following convention for the Fourier transform:
\[\hat f(\om)=\int_{\mbb R^d}\mr e^{\mr i \om\cdot x}f(x)\d x.\]
The Fourier transform of a fractional covariance kernel is given by
the corresponding spectral measure, which is defined as follows
(see, e.g., \cite[(1.10)]{ChenRoughSpace}, together with the fact
that $\tfrac{\Ga(2H+1)}{2\Ga(H)\Ga(1-H)}=\frac{\Ga(2H+1)\sin(H\pi)}{2\pi}$ for $H\in(0,1)$):
\begin{definition}
For every $H\in(0,1)$, define the measure
\begin{align}
\label{Equation: mu}
\dd\mu_H(x)=
\frac{\Ga(2H+1)}{2\Ga(H)\Ga(1-H)}|x|^{1-2H}\d x
\qquad x\in\mbb R.
\end{align}
Then, define $\mu=\mu_{H_1}\otimes\cdots\otimes\mu_{H_d}.$
For the time covariance, we denote
\begin{align}
\label{Equation: nu}
\nu=
\begin{cases}
\mu_{H_0}&\text{if }H_0<1,\\
\de_0&\text{if }H_0=1.
\end{cases}
\end{align}
\end{definition}

In particular, we have the Fourier representation
\[(\msf h\otimes\msf g)(s-t,x-y)
=\int_{\mbb R\times\mbb R^d}\mr e^{\mr i\ze(s-t)+\mr i\om\cdot(x-y)}\d(\nu\otimes\mu)(\ze,\om);\]
combining this with the well-known Fourier transform of the mollifiers
$f_\eps$ in Definition \eqref{Equation: Mollifier} (i.e., the Cauchy and Gaussian densities), we also have that
\[(\msf h\otimes\msf g)_\eps(s-t,x-y)
=\int_{\mbb R\times\mbb R^d}\mr e^{-\eps|\ze|-\eps|\om|^2/2}\mr e^{\mr i\ze(s-t)+\mr i\om\cdot(x-y)}\d(\nu\otimes\mu)(\ze,\om)\]
for $\eps>0$.
The purpose of this framework is to obtain the following convenient representations:
By Definition \ref{Definition: Approximate SILT and MILT}, one has
\begin{align}
\label{Equation: Beta eps Fourier}
\beta^{\eps}_t(B)=\int_{[0,t]^2}\int_{\mbb R\times\mbb R^d}\mr e^{-\eps(|\ze|+|\om|^2/2)}\mr e^{\mr i(s_2-s_1)\ze+\mr i\om\cdot(B(s_2)-B(s_1))}\d(\nu\otimes\mu)(\ze,\om)\dd s
\end{align}
and
\begin{align}
\label{Equation: Alpha eps Fourier}
\al^{\eps}_t(x)=\int_{[0,t]^2}\int_{\mbb R\times\mbb R^d}\mr e^{-\eps(|\ze|+|\om|^2/2)}\mr e^{\mr i(s_2-s_1)\ze+\mr i\om\cdot(B^{x_2}_2(s_2)-B^{x_1}_1(s_1))}\d(\nu\otimes\mu)(\ze,\om)\dd s.
\end{align}

\subsection{Proof of Proposition \ref{Proposition: Finite Epsilon Feynman-Kac}}
\label{Section: Proof of Proposition: Finite Epsilon Feynman-Kac}

By applying Tonelli's theorem in \eqref{Equation: Feynman-Kac}, we get
\[\msf M_{\ka,\eps_1}(t)\msf M_{\ka,\eps_2}(t)=\int_{D^2}\mbf E_B\left[\mbf 1_{\cap_{i\leq2}\{\tau_D(B_i^{x_i})>t\}}\mr e^{\sum_{i=1}^2\ka\int_0^t\xi_{\eps_i}(t-u,B_i^{x_i}(u))\d u}\right]\d x\]
for every $\eps_1,\eps_2,t>0$, where $\xi$, $B_1$, and $B_2$ are independent.
Then, by another application of Tonelli's theorem, we get the moment formulas
\begin{align}
\label{Equation: Finite Epsilon Feynman-Kac 1}
\mbf E\big[\msf M_{\ka,\eps}(t)\big]=\int_D\mbf E_B\left[\mbf 1_{\{\tau_D(B^x)>t\}}\mbf E_\xi\left[\mr e^{\ka\int_0^t\xi_\eps(t-u,B^x(u))\d u}\right]\right]\d x
\end{align}
and
\begin{align}
\label{Equation: Finite Epsilon Feynman-Kac 2}
\mbf E\big[\msf M_{\ka,\eps_1}(t)\msf M_{\ka,\eps_2}(t)\big]=\int_{D^2}\mbf E_B\left[\mbf 1_{\cap_{i\leq2}\{\tau_D(B_i^{x_i})>t\}}\mbf E_\xi\left[\mr e^{\sum_{i=1}^2\ka\int_0^t\xi_{\eps_i}(t-u,B_i^{x_i}(u))\d u}\right]\right]\d x,
\end{align}
where $\mbf E_\xi$ now denotes the expectation with respect to
$\xi$ only. Thus, the proof of Proposition \ref{Proposition: Finite Epsilon Feynman-Kac} is reduced to a moment generating function calculation.

Toward this end, we note that conditional on $B_i$, the random variables
\[\ka\int_0^t\xi_{\eps_i}\big(t-u,B^{x_i}_i(u)\big)\d u\] are jointly
Gaussian with mean zero. Moreover, by Fubini's theorem (by continuity, hence boundedness on compact sets of $u\mapsto\xi_\eps\big(t-u,B^x(u)\big)$), one has
\begin{multline*}
\mbf E_\xi\left[\left(\ka\int_0^t\xi_{\eps}\big(t-u,B^x(u)\big)\d u\right)^2\right]\\
=\ka^2\int_{[0,t]^2}\mbf E_\xi\big[\xi_\eps\big(t-u,B^x(u)\big)\xi_\eps\big(t-v,B^x(v)\big)\big]\d u\dd v
=\ka^2\be^\eps_t(B),
\end{multline*}
where the last equality follows from
\eqref{Equation: Covariance for eps1 and eps2} and Definition \ref{Definition: Approximate SILT and MILT}. Similarly,
\begin{multline*}
\mbf E_\xi\left[\left(\sum_{i=1}^2\ka\int_0^t\xi_{\eps_i}\big(t-u,B^{x_i}_i(u)\big)\d u\right)^2\right]\\
=\sum_{i=1}^2\ka^2\int_{[0,t]^2}\mbf E_\xi\big[\xi_{\eps_i}\big(t-u,B^{x_i}_i(u)\big)\xi_{\eps_i}\big(t-v,B^{x_i}_i(v)\big)\big]\d u\dd v\\
+2\ka^2\int_{[0,t]^2}\mbf E_\xi\big[\xi_{\eps_1}\big(t-u,B^{x_1}_1(u)\big)\xi_{\eps_2}\big(t-v,B^{x_2}_2(v)\big)\big]\d u\dd v,
\end{multline*}
which, by \eqref{Equation: Covariance for eps1 and eps2} and Definition \ref{Definition: Approximate SILT and MILT},
is equal to $\ka^2\sum_{i=1}^2\be^{\eps_i}_t(B_i)+2\ka^2\al_t^{(\eps_1+\eps_2)/2}(x)$.
Thus, we immediately obtain \eqref{Equation: Finite Epsilon Feynman-Kac E} from
\eqref{Equation: Finite Epsilon Feynman-Kac 1} thanks to the Gaussian moment generating function,
and we similarly obtain the following from \eqref{Equation: Finite Epsilon Feynman-Kac 2}:
\begin{align}
\label{Equation: Finite Epsilon Feynman-Kac 3}
\mbf E\big[\msf M_{\ka,\eps_1}(t)\msf M_{\ka,\eps_2}(t)\big]=\int_{D^2}\mbf E\left[\mbf 1_{\cap_{i\leq2}\{\tau_D(B_i^{x_i})>t\}}\mr e^{\frac{\ka^2}2\sum_{i=1}^2\be^{\eps_i}_t(B_i)+\ka^2\al_t^{(\eps_1+\eps_2)/2}(x)}\right]\d x.
\end{align}
With this in hand, we obtain \eqref{Equation: Finite Epsilon Feynman-Kac Var} by combining
\eqref{Equation: Finite Epsilon Feynman-Kac 3} (in the special case $\eps_1=\eps_2=\eps$),
\eqref{Equation: Finite Epsilon Feynman-Kac E}, and $\mbf{Var}\big[\msf M_{\ka,\eps}(t)\big]
=\mbf E\big[\msf M_{\ka,\eps}(t)^2\big]-\mbf E\big[\msf M_{\ka,\eps}(t)\big]^2$.

\subsection{References/Proofs for Propositions \ref{Proposition: Limit and UI of SILT}
and \ref{Proposition: Limit and UI of MILT}}
\label{Section: Proof of Limit and UI of SILT and MILT}

We begin with the proof of Proposition \ref{Proposition: Limit and UI of MILT},
since we use some of it in the proof of Proposition \ref{Proposition: Limit and UI of SILT}.

\subsubsection{Proposition \ref{Proposition: Limit and UI of MILT}}

We follow the strategy in \cite[Section 2]{ChenRoughSpace},
which proved the result in the special case $x=0$.
Let $t>0$ be fixed.
For every $x\in(\mbb R^d)^2$, $m\in\mbb N$, $\eps\in(0,\infty)^m$, $\ze\in\mbb R^m$, and $\om\in(\mbb R^d)^m$,
let us denote
\begin{align}
\label{Equation: A Moment Function}
\msf A^{(m)}_\eps(\ze,\om)=\mr e^{-\sum_{k=1}^m\eps_k(|\ze_k|+|\om_k|^2/2)}\end{align}
and
\begin{multline}
\label{Equation: B Moment Function}
\msf B^{(m)}_t(\ze,\om)
=\int_{[0,t]^{2m}}\mr e^{\mr i\sum_{k=1}^m(s_k-v_k)\ze_k}\mbf E\left[\mr e^{\mr i\sum_{k=1}^m\om_k\cdot(B_2(s_k)-B_1(v_k))}\right]\d s\dd v\\
=\left|\int_{[0,t]^{m}}\mr e^{\mr i\sum_{k=1}^ms_k\ze_k}\mbf E\left[\mr e^{\mr i\sum_{k=1}^m\om_k\cdot B(s_k)}\right]\d s\right|^2,
\end{multline}
where the last equality follows from the independence of $B_1$ and $B_2$.
Then, a straightforward application of Tonelli's theorem in \eqref{Equation: Alpha eps Fourier}
yields
\[\mbf E\big[\al^{\eps_1}_t(x)\cdots\al^{\eps_m}_t(x)\big]
=\int_{\mbb R^m\times(\mbb R^d)^m}\mr e^{\mr i(x_2-x_1)\cdot\sum_{k=1}^m\om_k}
\msf A^{(m)}_\eps(\ze,\om)\msf B^{(m)}_t(\ze,\om)\d\nu^{\otimes m}(\ze)\dd\mu^{\otimes m}(\om).\]
If we show that the above integral converges to a finite limit as $\eps_i\to0$
for $i=1,\ldots,m$ for every $m\geq1$, then it will follow that $\al^\eps_t(x)$ converges in $L^p$ for all $p\geq1$
to a random variable $\al_t(x)$.
For this, given that $\msf A^{(m)}_\eps(\ze,\om)\to1$ pointwise as $\eps\to0$ and that
\begin{align}
\label{Equation: MILT Proof 1}
\left|\msf A^{(m)}_\eps(\ze,\om)\mr e^{\mr i(x_2-x_1)\cdot\sum_{k=1}^m\om_k}\right|\leq1,
\end{align}
it suffices to show by the dominated convergence theorem that
the integral of $\msf B$ is finite. By \cite[(3.1) and (4.1)]{ChenRoughSpace},
under Assumption \ref{Assumption},
there exists a constant $C>0$ independent of $t$ such that
\begin{align}
\label{Equation: Moment Bouns for MILT}
\int_{\mbb R^m\times(\mbb R^d)^m}\msf B^{(m)}_t(\ze,\om)\d\nu^{\otimes m}(\ze)\dd\mu^{\otimes m}(\om)\leq C^m(m!)^{d-H}t^{(1+\rho)m},
\qquad m\in\mbb N.
\end{align}
Thus, \eqref{Equation: MILT Limit} holds. This bound also 
implies \eqref{Equation: MILT UI}:
Since $\al^\eps_t(x)$ is nonnegative (by Definition \ref{Definition: Approximate SILT and MILT}
and the fact that $\msf h\otimes g$ is a positive measure under Assumption \ref{Assumption}),
for every $m\geq1$,
\begin{multline}
\mbf E\big[\al^{\eps}_t(x)^m\big]
=\left|\int_{\mbb R^m\times(\mbb R^d)^m}\mr e^{\mr i(x_2-x_1)\cdot\sum_{k=1}^m\om_k}
\msf A^{(m)}_\eps(\ze,\om)\msf B^{(m)}_t(\ze,\om)\d\nu^{\otimes m}(\ze)\dd\mu^{\otimes m}(\om)\right|\\
\leq\int_{\mbb R^m\times(\mbb R^d)^m}\msf B^{(m)}_t(\ze,\om)\d\nu^{\otimes m}(\ze)\dd\mu^{\otimes m}(\om),
\end{multline}
where the inequality follows from the triangle inequality (i.e., taking an absolute value inside an integral) and \eqref{Equation: MILT Proof 1}.
Combining this with \eqref{Equation: Moment Bouns for MILT}, we are led to
\[\sup_{\eps>0,~x\in(\mbb R^d)^2}\mbf E\Big[\exp\big(c\al^\eps_t(x)\big)\Big]
\leq\sum_{m=0}^\infty\frac{(cCt^{1+\rho})^m}{(m!)^{1-(d-H)}};\]
this is finite for all $t>0$ when $d-H<1$, and finite for small enough $t>0$
when $d-H=1$, as desired. The corresponding bound for $\al_t(x)$'s exponential
moments then follows from \eqref{Equation: MILT Limit} and uniform integrability.
The proof of Proposition \ref{Proposition: Limit and UI of MILT} is thus complete.

\subsubsection{Proposition \ref{Proposition: Limit and UI of SILT}}

We use a standard triangle decomposition argument (e.g.,
\cite[Figure 1]{ChenBook}), together
with estimates that appeared in \cite[Lemma 4.4]{ChenDeyaOuyangTindel}:
Let $t>0$ be fixed throughout this proof.
For every integer $N\geq1$
and $1\leq k\leq 2^{N-1}$, we let
\begin{align}
\label{Equation: J and I Intervals}
J_{N,k}=\left[\frac{(2k-2)t}{2^N},\frac{(2k-1)t}{2^N}\right)
\qquad\text{and}\qquad
I_{N,k}=\left[\frac{(2k-1)t}{2^N},\frac{2kt}{2^N}\right).
\end{align}
Given that the products $J_{N,k}\times I_{N,k}$ form a disjoint union (up to a Lebesgue null set) of
the set
$\{(s_1,s_2)\in[0,t]:s_1<s_2\}$, we can write
\begin{align}
\label{Equation: beta triangle decomposition}
\be^\eps_t(B)=2\sum_{N=1}^\infty\sum_{k=1}^{2^{N-1}}a^\eps_{N,k},
\end{align}
where
\[a_{N,k}^\eps=\int_{J_{N,k}\times I_{N,k}}(\msf h\otimes\msf g)_\eps\big(s_2-s_1,B(s_2)-B(s_1)\big)\d s.\]
If we consider the midpoint $(2k-1)t/2^N$ of $J_{N,k}\cup I_{N,k}$, and run the Brownian motion forward and backward in time, then we get
the equality in distribution
\[a_{N,k}^\eps\deq\int_{[0,t/2^N)^2}(\msf h\otimes\msf g)_\eps\big(s_2+s_1,B_2(s_2)-B_1(s_1)\big)\d s.\]
Moreover, for every $N\ge1$,
the variables $a^\eps_{N,1},\ldots,a^\eps_{N,2^{N-1}}$ are independent because they involve disjoint Brownian increments.

Now,
if we define
\begin{multline*}
\tilde{\msf B}^{(m)}_t(\ze,\om)
=\int_{[0,t]^{2m}}\mr e^{\mr i\sum_{k=1}^m(s_k+v_k)\ze_k}\mbf E\left[\mr e^{\mr i\sum_{k=1}^m\om_k\cdot(B_2(s_k)-B_1(v_k))}\right]\d s\dd v\\
=\left(\int_{[0,t]^{m}}\mr e^{\mr i\sum_{k=1}^ms_k\ze_k}\mbf E\left[\mr e^{\mr i\sum_{k=1}^m\om_k\cdot B(s_k)}\right]\d s\right)^2,
\end{multline*}
then the same Fourier transform argument
used in the previous section implies that
\begin{align}
\label{Equation: Triangle Decomp. Moment Bound}
\mbf E\big[a_{N,k}^{\eps_1}\cdots a_{N,k}^{\eps_m}\big]
=\int_{\mbb R^m\times(\mbb R^d)^m}
\msf A^{(m)}_\eps(\ze,\om)\tilde{\msf B}^{(m)}_{t/2^N}(\ze,\om)\d\nu^{\otimes m}(\ze)\dd\mu^{\otimes m}(\om)
\end{align}
(recall \eqref{Equation: A Moment Function} and \eqref{Equation: B Moment Function}). Given that this moment is nonnegative (since $a_{N,k}^\eps$ is itself nonnegative) and $|\tilde{\msf B}^{(m)}_t(\ze,\om)|=\msf B^{(m)}_t(\ze,\om)$, it follows from
the triangle inequality and \eqref{Equation: Moment Bouns for MILT} that for every $N\geq1$,
\begin{align}
\label{Equation: aNk moment bounds}
\mbf E\big[a_{N,k}^{\eps_1}\cdots a_{N,k}^{\eps_m}\big]\leq\mbf E[\al_{t/2^N}(0)^m]
\leq C^m(m!)^{d-H}(t/2^N)^{(1+\rho)m}
\end{align}
for a constant $C>0$ independent of $t$ and $N$. In particular, the same
argument we have used for
Proposition \ref{Proposition: Limit and UI of MILT} shows that for every choice of $N$ and $k$, there exists a random variable $a_{N,k}$ that is the $L^p$ limit of $a_{N,k}^{\eps}$ as $\eps\to0$ for all $p\geq1$.

We now move on to the construction of the limits $\ga_t(B)$ in the statement of Proposition \ref{Proposition: Limit and UI of SILT}, and
the proofs of the corresponding exponential moment bounds in
\eqref{Equation: subcritical SILT UI} and \eqref{Equation: renormalized SILT UI}.
First consider the case where $\rho>0$. Note that by Fubini's theorem
(which we can apply thanks to the exponential weight $\mr e^{-\eps(|\ze|+|\om|^2/2)}$),
we can write \eqref{Equation: Beta eps Fourier} as
\begin{align}
\label{Equation: Beta eps Fourier conj}
\beta^{\eps}_t(B)=\int_{\mbb R\times\mbb R^d}\mr e^{-\eps(|\ze|+|\om|^2/2)}\left|\int_0^t\mr e^{\mr is\ze+\mr i\om\cdot B(s)}\d s\right|^2\d(\nu\otimes\mu)(\ze,\om).
\end{align}
These variables are monotonically increasing as $\eps\downarrow0$.
Thus, by the monotone convergence theorem, we can define
\[\ga_t(B)=\int_{\mbb R\times\mbb R^d}\left|\int_0^t\mr e^{\mr is\ze+\mr i\om\cdot B(s)}\d s\right|^2\d(\nu\otimes\mu)(\ze,\om)\]
and have that $\be^\eps_t(B)\to\ga_t(B)$ surely,
although, a priori, this limiting random element could be infinite with positive probability.
By the Beppo Levi theorem (and monotonicity),
\[\mbf E\big[\ga_t(B)^m\big]=\sup_{\eps>0}\mbf E\big[\be^\eps_t(B)^m\big]\qquad\text{for every }m\geq1.\]
Therefore, by \eqref{Equation: beta triangle decomposition}, \eqref{Equation: Triangle Decomp. Moment Bound},
and Minkowski's inequality, we have
\begin{multline*}
\mbf E[\ga_t(B)^m]^{1/m}
\leq\sup_{\eps>0}2\sum_{N=1}^\infty\sum_{k=1}^{2^{N-1}}\mbf E[(a^\eps_{N,k})^m]^{1/m}\\
\leq 2C(m!)^{(d-H)/m}t^{1+\rho}\sum_{N=1}^\infty2^{N-1}\cdot 2^{-N(1+\rho)}
=\frac{C(m!)^{(d-H)/m}t^{1+\rho}}{(2^\rho-1)},
\end{multline*}
where we can apply the geometric series to $2^{-N\rho}$ because $\rho>0$.
Thus, $\ga_t(B)$ is almost surely finite.
Moreover, if we let $\tilde C=\frac{Ct^{1+\rho}}{(2^\rho-1)}$, then this moment estimate implies that
\[\sup_{\eps>0}\mbf E\left[\exp\big(c\be^\eps_t(B)\big)\right]\leq\mbf E\left[\exp\big(c\ga_t(B)\big)\right]
=\sum_{m=0}^\infty\frac{c^m\mbf E[\ga_t(B)^m]}{m!}
\leq\sum_{m=0}^\infty(\tilde Cc)^m(m!)^{(d-H)-1}\]
for every $c>0$. When $\rho>0$, it is necessarily
the case that $d-H<1$ (since $2H_0-1\leq1$).
In particular, for every $\eta>0$, there exists a constant $C_\eta>0$ large enough so that
\begin{align}
\label{Equation: (d-H)<1 Factorial Bound}
(m!)^{d-H}\leq C_\eta\eta^mm!,\qquad m\geq1.
\end{align}
Taking $\eta$
small enough so that $\tilde Cc\eta<1$, this yields
\[\sup_{\eps>0}\mbf E\left[\exp\big(c\be^\eps_t(B)\big)\right]\leq\mbf E\left[\exp\big(c\ga_t(B)\big)\right]
\leq C_\eta\sum_{m=0}^\infty(\tilde Cc\eta)^m
=\frac{C_\eta}{1-\tilde Cc\eta}<\infty,\]
as desired.

Consider now the case $\rho\leq0$.
Our first aim is to show that
\begin{align}
\label{Equation: Renormalized Gamma Triangle Decomposition}
\ga_t(B)=
2\sum_{N=1}^\infty\sum_{k=1}^{2^{N-1}}\big(a_{N,k}-\mbf E[a_{N,k}]\big)
\end{align}
satisfies the desired $L^2$ limit.
For this purpose, let us define the quantities
\begin{align}
\label{Equation: Gamma Epsilon}
\ga^\eps_t(B)=
\begin{cases}
\be^\eps_t(B)&\rho>0,\\
\be^\eps_t(B)-\mbf E\big[\be^\eps_t(B)\big]&\rho\leq0,
\end{cases}
\qquad\eps>0,
\end{align}
and, assuming $\rho\leq0$,
\[\ga^\eps_{M;t}(B)=
2\sum_{N=1}^M\sum_{k=1}^{2^{N-1}}\big(a^\eps_{N,k}-\mbf E[a^\eps_{N,k}]\big),\qquad\eps>0,~M\geq1.\]
Define $\ga_{M;t}(B)$ in the same way, but with every instance of $a^\eps_{N,k}$ replaced by the corresponding limit $a_{N,k}$.
Then, by the triangle inequality,
\begin{multline}
\label{Equation: From Triangle Convergence to Full Convergence}
\mbf E\big[\big|\ga^\eps_t(B)-\ga_t(B)\big|^2\big]^{1/2}
\leq
\mbf E\big[\big|\ga^\eps_t(B)-\ga^\eps_{M;t}(B)\big|^2\big]^{1/2}+\mbf E\big[\big|\ga_{M;t}(B)-\ga_t(B)\big|^2\big]^{1/2}\\
+\mbf E\big[\big|\ga^\eps_{M;t}(B)-\ga_{M;t}(B)\big|^2\big]^{1/2}.
\end{multline}
For every $M\geq1$,
the contribution on the second line of
\eqref{Equation: From Triangle Convergence to Full Convergence} vanishes
 as $\eps\to0$ by the $L^2$ convergence
 of $a_{N,k}^\eps$ to $a_{N,k}$ for each fixed $N$ and $k$. Thus, it suffices to show that the two quantities on the first
 line of the right-hand side of \eqref{Equation: From Triangle Convergence to Full Convergence} can be made negligible as $M\to\infty$, uniformly in $\eps>0$.
 Note that for every $M\geq1$,
\[\ga^\eps_{t}(B)-\ga^\eps_{M;t}(B)=
2\sum_{N=M+1}^\infty\sum_{k=1}^{2^{N-1}}\big(a^\eps_{N,k}-\mbf E[a^\eps_{N,k}]\big),\]
and similarly for $\ga_{t}(B)-\ga_{M;t}(B)$. Thus, by a combination of the triangle inequality and the fact that the $a_{N,k}$'s are i.i.d. for fixed $N$ as $k$ varies, it suffices to prove that
\[\sup_{\eps>0}
\sum_{N=1}^\infty\sqrt{2^{N-1}\mbf{Var}[a^\eps_{N,1}]}<\infty
\qquad\text{and}\qquad
\sum_{N=1}^\infty\sqrt{2^{N-1}\mbf{Var}[a_{N,1}]}<\infty.\]
By combining the fact that $\mbf{Var}[X]\leq\mbf E[X^2]$
and the moment bound \eqref{Equation: aNk moment bounds} (which also holds for $a_{N,1}$ by virtue of being the $L^2$ limit of $a^\eps_{N,1}$), for this it suffices that
\[\sum_{N=1}^\infty 2^{(N-1)/2}2^{-(1+\rho)N}=2^{-1/2}\sum_{N=1}^\infty2^{-(1/2+\rho)N}<\infty.\]
This holds even when $\rho\leq0$ because \eqref{Equation: Assumption} implies that $\rho>-1/2$.

It now only remains to prove \eqref{Equation: renormalized SILT UI}. Following the argument
in \cite[paragraph containing (4.25)]{ChenDeyaOuyangTindel}
(which is a standard H\"older inequality trick
using the triangle decomposition; e.g. \cite[Proof of Theorem 2.4.2]{ChenBook} or \cite[Step 2 of the Proof of Theorem 4.6]{HuHuangNualartTindel}),
we get the following dichotomy:
\begin{enumerate}
\item If, for every $\eta>0$, there exists
a constant $C_\eta>0$ large enough so that
\[\sup_{\eps>0}\mbf E[(a^\eps_{N,k})^m]\leq C_\eta(C\eta t/2^N)^{qm}m!,\qquad m\geq1,\]
for some $C>0$ and $q>1/2$,
then \eqref{Equation: renormalized SILT UI} holds for all $t>0$.
\item If there exists $C>0$ and $q>1/2$ such that
\[\sup_{\eps>0}\mbf E[(a^\eps_{N,k})^m]\leq (Ct/2^N)^{qm}m!,\qquad m\geq1,\]
then \eqref{Equation: renormalized SILT UI} holds for all $t$ small enough.
\end{enumerate}
Given that $\rho>-1/2$ (and thus $q=(1+\rho)>1/2$), item (1)
holds whenever $d-H<1$ by
\eqref{Equation: aNk moment bounds} and
\eqref{Equation: (d-H)<1 Factorial Bound},
and item (2) holds when $d-H=1$ by \eqref{Equation: aNk moment bounds}.

\subsection{Proof of Lemma \ref{Lemma: Scaling of Beta and Alpha}}
\label{Section: Proof of Scaling of Beta and Alpha}

By Brownian scaling, for every $x\in\mbb R^d$,
\[B^x(\cdot)\deq\sqrt tB^{x/\sqrt t}(\cdot/t).\]
Thus, by \eqref{Equation: Beta eps Fourier},
\[\be^\eps_t(B)\deq\int_{[0,t]^2}\int_{\mbb R\times\mbb R^d}\mr e^{-\eps(|\ze|+|\om|^2/2)}\mr e^{\mr i(s_2-s_1)\ze+\mr i(\sqrt t\om)\cdot(B(s_2/t)-B(s_1/t))}\d(\nu\otimes\mu)(\ze,\om)\dd s.\]
By the change of variables $v=s/t$, $\dd s=t^2\dd v$, we get
\[\be^\eps_t(B)\deq t^2\int_{[0,1]^2}\int_{\mbb R\times\mbb R^d}\mr e^{-\eps(|\ze|+|\om|^2/2)}\mr e^{\mr i(v_2-v_1)t\ze+\mr i(\sqrt t\om)\cdot(B(v_2)-B(v_1))}\d(\nu\otimes\mu)(\ze,\om)\dd v.\]
Then, by the change of variables $(\chi,\ga)=(t\ze,\sqrt t\om)$,
\[\d(\nu\otimes\mu)(\ze,\om)=t^{-2+2H_0-(d-H)}\d(\nu\otimes\mu)(\chi,\ga)\]
(the latter of which follows from \eqref{Equation: mu} and \eqref{Equation: nu}),
we get
\begin{multline*}
\be^\eps_t(B)
\deq t^{1+(2H_0-1)-(d-H)}\\
\times\int_{[0,1]^2}\int_{\mbb R\times\mbb R^d}\mr e^{-(\eps/t)(|\chi|+|\ga|^2/2)}\mr e^{\mr i(v_2-v_1)\chi+\mr i\ga\cdot(B(v_2)-B(v_1))}\d(\nu\otimes\mu)(\chi,\ga)\dd v
=t^{1+\rho}\be_1^{\eps/t}(B).
\end{multline*}
The proof of the scaling property for $\al^\eps_t(x)$ is identical,
except that the additional $x$-dependence in the
term of the form $\mr e^{\mr i\om\cdot(B^{x_2}_2(s_2)-B^{x_1}_1(s_1))}$ in \eqref{Equation: Alpha eps Fourier}
gets scaled into $\mr e^{\mr i\ga\cdot(B^{x_2/\sqrt t}_2(u_2)-B^{x_1/\sqrt t}_1(u_1))}$.

\subsection{Proof of Proposition \ref{Proposition: SILT Expectation Limits}}
\label{Section: SILT Expectation Limits}

\subsubsection{Outline of Proof}

Since the proof of Proposition \ref{Proposition: SILT Expectation Limits} involves many different
steps, we begin by stating several technical lemmas without proof, and then use the lemmas in question
to get the desired conclusion. Then, we prove the technical lemmas one at a time.

Our first technical lemma consists of deriving integral formulas for $\mbf E\big[\tfrac{\ka^2}{2}\be^\th_1(B)\big]$ that are amenable
to computation, namely:

\begin{lemma}[Moment Formulas]
\label{Lemma: Moment Formulas}
Let $\ka,\th>0$.
Let us denote
\begin{align}
\label{Equation: J function}
\msf J(\th)=
\begin{cases}
\displaystyle\int_{-1}^1(1-|v|)\int_{\mbb R^2}\mr e^{-\th(|a|+|b|)}\mr e^{\mr iva-|v|\cdot|b|}|a|^{1-2H_0}|b|^{d-H-1}\d a\dd b\dd v,&H_0<1,
\vspace{5pt}\\
\displaystyle\int_{-1}^1(1-|v|)\int_{\mbb R}\mr e^{-\th|b|}\mr e^{-|v|\cdot|b|}|b|^{d-H-1}\d b\dd v,&H_0=1,
\end{cases}
\end{align}
and the constant
\begin{align}
\label{Equation: a constant}
\mf a=
\begin{cases}
\displaystyle\frac{1}{2^{H+3}\Ga(d-H)\Ga(1-H_0)}\prod_{j=0}^d\frac{\Ga(2H_j+1)}{\Ga(H_j)}&H_0<1,
\vspace{5pt}\\
\displaystyle\frac{1}{2^{H+2}\Ga(d-H)}\prod_{j=1}^d\frac{\Ga(2H_j+1)}{\Ga(H_j)}&H_0=1.
\end{cases}
\end{align}
In all cases,
\begin{align}
\label{Equation: Moment Formula}
\mbf E\big[\tfrac{\ka^2}{2}\be^\th_1(B)\big]=\ka^2\mf a\,\msf J(\th).
\end{align}
\end{lemma}

Lemma \ref{Lemma: Moment Formulas} is proved in Section \ref{Section: Proof of Lemma: Moment Formulas},
and follows from a simplification of the expectation of the Fourier integral \eqref{Equation: Beta eps Fourier} using
various changes of variables.
Next, we must analyze the integral $\msf J(\th)$. The behavior of the latter depends significantly on
$\rho$, and so we split the results into three parts:

\begin{lemma}[Integral Asymptotics 1]
\label{Lemma: Integral Asymptotics - Subcritical}
When $\rho>0$ (which necessarily
implies that $H_0>1/2$), one has
\begin{align}
\label{Equation: Integral Asymptotics - Subcritical}
\lim_{\th\to0}\msf J(\th)
=\begin{cases}
\displaystyle\frac{4\Ga(H_0)\Ga(1-H_0)\Ga(d-H)}{\Ga(2H_0-1)(1+\rho)\rho}&H_0<1,
\vspace{5pt}\\
\displaystyle\frac{4\Ga(d-H)}{(1+\rho)\rho}&H_0=1.
\end{cases}
\end{align}
\end{lemma}

\begin{lemma}[Integral Asymptotics 2]
\label{Lemma: Integral Asymptotics - Critical}
If $\rho=0$ (which also implies $H_0>1/2$),
there exists a continuous and bounded function $\msf b:(0,\infty)\to\mbb R$ such that
\begin{align}
\label{Equation: Integral Asymptotics - Critical}
\msf J(\th)=
\begin{cases}
4\Ga(H_0)\Ga(1-H_0)\ell(\th)+\msf b(\th)&1/2<H_0<1,\\
4\big((1+\th)\log(1+1/\th)-1\big)& H_0=1,
\end{cases}
\end{align}
where we recall the definition of $\ell$ in \eqref{Equation: ell function}.
\end{lemma}

\begin{lemma}[Integral Asymptotics 3]
\label{Lemma: Integral Asymptotics - Supercritical}
Suppose that $\rho<0$, which, under Assumption \ref{Assumption}, implies that $H_0<1$. 
There is a bounded continuous function $\msf z_0:(0,\infty)\to\mbb R$
such that
for every $\th>0$, one has
\begin{align}
\label{Equation: Integral Asymptotics - Supercritical}
\msf J(\th)=8\Ga(-\rho)\left(\int_0^{\pi/2}\frac{(\cos\theta)^{-(2H_0-1)}(\sin\theta)^{d-H}}{\big(\cos\theta+\sin\theta\big)^{-\rho}}\d\theta\right)\th^{\rho}+\msf z_0(\th).
\end{align}
Moreover, one has
\begin{align}
\label{Equation: Integral Asymptotics - Supercritical Zero Limit}
\lim_{\th\to0}\msf z_0(\th)=
\begin{cases}
\displaystyle \frac{4\Ga(H_0)\Ga(1-H_0)\Gamma(d-H)}{\Ga(2H_0-1)(1+\rho)\rho}&H_0>1/2,\\
0&H_0=1/2.
\end{cases}
\end{align}
\end{lemma}

With this in hand, we may now prove
Proposition \ref{Proposition: SILT Expectation Limits}:
Firstly, we obtain \eqref{Equation: Separate Expectation Limits - Subcritical} by combining \eqref{Equation: Moment Formula},
\eqref{Equation: Integral Asymptotics - Subcritical}, and the fact that
\begin{align}
\label{Equation: a to b}
\frac{\ka^2\mf b}{(1+\rho)\rho}=\begin{cases}
\displaystyle\ka^2\mf a\frac{4\Ga(H_0)\Ga(1-H_0)\Ga(d-H)}{\Ga(2H_0-1)(1+\rho)\rho}&H_0<1,
\vspace{5pt}\\
\displaystyle\ka^2\mf a\frac{4\Ga(d-H)}{(1+\rho)\rho}&H_0=1,
\end{cases}
\end{align}
by definitions of $\mf a$ and $\mf b$ in
\eqref{Equation: a constant} and \eqref{Equation: b Constant},
the simplification $\frac{\Gamma(2H_0+1)}{\Gamma(2H_0-1)}=2H_0(2H_0-1)$ when $H_0<1$,
and the fact that $H_0(2H_0-1)=1$ when $H_0=1$.
Secondly, we obtain \eqref{Equation: Separate Expectation Limits - Critical} by combining \eqref{Equation: Moment Formula}, \eqref{Equation: Integral Asymptotics - Critical},
and \eqref{Equation: a to b}, noting that when $\rho=0$,
$\Gamma(d-H)=\Gamma(2H_0-1)$, and if we additionally assume
that $H_0=1$, then $\Gamma(d-H)=\Gamma(1)=1$.
In particular, the functions $\msf b_i$ are defined as
\[\msf b_1(\th)=\mf a\msf b(\th)
\qquad\text{and}\qquad
\msf b_2(\th)=-4\mf a.\]
Thirdly, we get \eqref{Equation: Separate Expectation Limits - Supercritical} by combining
\eqref{Equation: Moment Formula}, \eqref{Equation: Integral Asymptotics - Supercritical},
\eqref{Equation: Integral Asymptotics - Supercritical Zero Limit}, and \eqref{Equation: a to b} (in the case $H_0<1$),
together with
\[\ka^2\mf a8\Ga(-\rho)\left(\int_0^{\pi/2}\frac{(\cos\theta)^{-(2H_0-1)}(\sin\theta)^{d-H}}{\big(\cos\theta+\sin\theta\big)^{-\rho}}\d\theta\right)=\ka^2\mf d\]
by definitions of $\mf a$ and $\mf d$ in \eqref{Equation: a constant} and \eqref{Equation: Renormalization - Supercritical constant}.

At this point, in order to complete the proof
of Proposition \ref{Proposition: SILT Expectation Limits},
it only remains to prove Lemmas \ref{Lemma: Moment Formulas}--\ref{Lemma: Integral Asymptotics - Supercritical}.
We now turn to this task.

\subsubsection{Proof of Lemma \ref{Lemma: Moment Formulas}}
\label{Section: Proof of Lemma: Moment Formulas}

By Tonelli's theorem and \eqref{Equation: Beta eps Fourier},
\[\mbf E\big[\tfrac{\ka^2}2\beta^{\th}_1(B)\big]=\frac{\ka^2}{2}\int_{[0,1]^2}\int_{\mbb R\times\mbb R^d}\mr e^{-\th(|\ze|+|\om|^2/2)}\mr e^{\mr i(s_2-s_1)\ze}\mbf E\left[\mr e^{\mr i\om\cdot(B(s_2)-B(s_1))}\right]\d(\nu\otimes\mu)(\ze,\om)\dd s.\]
Given that $B(s_2)-B(s_1)\sim N\big(0,|s_2-s_1|I_d\big)$,
where $N(m,\Si)$ denotes the multivariate Gaussian distribution with mean vector $m$
and covariance matrix $\Sigma$, and $I_d$ denotes the $d\times d$ identity matrix,
this becomes
\begin{align}
\label{Equation: Moment Formulas 1}
\mbf E\big[\tfrac{\ka^2}2\beta^{\th}_1(B)\big]=\frac{\ka^2}{2}\int_{[0,1]^2}\int_{\mbb R\times\mbb R^d}\mr e^{-\th(|\ze|+|\om|^2/2)}\mr e^{\mr i(s_2-s_1)\ze-|s_2-s_1|\cdot|\om|^2/2}\d(\nu\otimes\mu)(\ze,\om)\dd s.
\end{align}

If we isolate the $\dd\om$ integral above, expanding $\mu$ according to \eqref{Equation: mu},
then we get the contribution
\[\frac{1}{2^d}\prod_{j=1}^d\frac{\Ga(2H_j+1)}{\Ga(H_j)\Ga(1-H_j)}\int_{\mbb R^d}
\mr e^{-\th|\om|^2/2-|s_2-s_1|\cdot|\om|^2/2}\prod_{j=1}^{d}|\om_j|^{1-2H_j}\d\om.\]
By a spherical change of variables, this integral can be split into the product of two contributions:
The angular component
\begin{multline*}
\mf C_{\mr{ang}}=\int_{[0,\pi]^{d-2}\times[0,2\pi]}(\sin\theta_1)^{d-2}(\sin\theta_2)^{d-3}\cdots\sin\theta_{d-2}\\
\times|\cos\theta_1|^{1-2H_1}|\sin\theta_1\cos\theta_2|^{1-2H_2}\cdots|\sin\theta_1\cdots\sin\theta_{d-2}\cos\theta_{d-1}|^{1-2H_{d-1}}\cdots\\
\times|\sin\theta_1\cdots\sin\theta_{d-2}\sin\theta_{d-1}|^{1-2H_d}\d\theta_1\cdots\dd\theta_{d-1},
\end{multline*}
when $d\geq2$ or $\mf C_{\mr{ang}}=2$ when $d=1$,
and the radial component
\begin{align}
\label{Equation: Radial Component}
\frac{1}{2^d}\prod_{j=1}^d\frac{\Ga(2H_j+1)}{\Ga(H_j)\Ga(1-H_j)}\int_0^\infty
\mr e^{-\th r^2/2-|s_2-s_1| r^2/2} r^{2(d-H)-1}\d r.
\end{align}
Our aim is to calculate $\mf C_{\mr{ang}}$, simplify the radial component,
and then put these contributions back into \eqref{Equation: Moment Formulas 1}.

Regarding $\mf C_{\mr{ang}}$,
by the
well-known formula for the fractional absolute moments of the univariate Gaussian distribution,
i.e.,
\begin{align}
\label{Equation: Gaussian Fractional Moment}
\mbf E\big[|N(0,\si^2)|^p\big]=\tfrac{\si^p2^{p/2}\Ga((p+1)/2)}{\sqrt\pi},\qquad p\in(-1,\infty),
\end{align}
one has
\[\int_{\mbb R^d}\frac{\mr e^{-|\om|^2/2}}{2^{d-H}}\prod_{j=1}^d|\om_j|^{1-2H_j}\d\om=\prod_{j=1}^d\Ga (1-H_j).\]
If we compute this same integral but with a spherical change of variables first (together with
a second application of \eqref{Equation: Gaussian Fractional Moment} to calculate the radial part of the integral), then we also have that
\[\int_{\mbb R^d}\frac{\mr e^{-|\om|^2/2}}{2^{d-H}}\prod_{j=1}^d|\om_j|^{1-2H_j}\d\om=\mf C_{\mr{ang}}\int_0^\infty\frac{\mr e^{-r^2/2}}{2^{d-H}} r^{2(d-H)-1}\d r
=\mf C_{\mr{ang}}\frac{\Ga(d-H)}{2}.\]
We therefore conclude that
\begin{align}
\label{Equation: Cang Value}
\mf C_{\mr{ang}}=2\frac{\prod_{j=1}^d\Ga (1-H_j)}{\Ga(d-H)}.
\end{align}
Regarding the radial component,
by the change of variables $b=r^2/2$, $\dd b=r\dd r$,
and then artificially extending the domain of integration from $b\in[0,\infty)$
to $b\in\mbb R$ by symmetry, we have the identity
\[\eqref{Equation: Radial Component}
=
\frac{1}{2^{H+2}}\prod_{j=1}^d\frac{\Ga(2H_j+1)}{\Ga(H_j)\Ga(1-H_j)}\int_{\mbb R}
\mr e^{-\th|b|-|s_2-s_1|\cdot|b|}|b|^{d-H-1}\d b.\]
Putting everything back into 
\eqref{Equation: Moment Formulas 1}, we arrive at
\begin{multline*}
\mbf E\big[\tfrac{\ka^2}2\beta^{\th}_1(B)\big]=\frac{\ka^2}{2^{H+2}\Ga(d-H)}\prod_{j=1}^d\frac{\Ga(2H_j+1)}{\Ga(H_j)}\\
\times\int_{[0,1]^2}\int_{\mbb R^2}\mr e^{-\th(|\ze|+|b|)}\mr e^{\mr i(s_2-s_1)\ze-|s_2-s_1|\cdot|b|}|b|^{d-H-1}\d\nu(\ze)\dd b\dd s.
\end{multline*}
We then get \eqref{Equation: Moment Formula}
by \eqref{Equation: nu} (and the cosmetic change $a=\ze$),
and the change of variables $(u_1,v)=(u_1,u_2-u_1)$,
whereby
for any nonnegative or integrable measurable function $f$, one has
\begin{align}
\label{Equation: u2-u1 to Uniform Integral Trick}
\int_{[0,1]^2}f(u_2-u_1)\d u=\int_{-1}^1(1-|v|)f(v)\d v.
\end{align}

\subsubsection{Proof of Lemma \ref{Lemma: Integral Asymptotics - Subcritical}}
\label{Section: Proof of Lemma: Integral Asymptotics - Subcritical}

We begin with some generalities that hold for any $-1/2<\rho<1$.
For every fixed $\th>0$, the integrand in $\msf J(\th)$
is integrable thanks to the exponential factor $\mr e^{-\th(|a|+|b|)}$.
Thus, by Fubini's theorem,
\[\msf J(\th)
=\int_{-1}^1(1-|v|)\left(\int_{\mbb R}\mr e^{-\th|a|}\mr e^{\mr iva}|a|^{1-2H_0}\d a\right)
\left(\int_{\mbb R}\mr e^{-\th|b|}\mr e^{-|v|\cdot|b|}|b|^{d-H-1}\d b\right)\dd v\]
when $H_0<1$, and
\begin{align*}
\msf J(\th)
=\int_{-1}^1(1-|v|)\int_{\mbb R}\mr e^{-\th|b|}\mr e^{-|v|\cdot|b|}|b|^{d-H-1}\d b\dd v
\end{align*}
when $H_0=1$.
On the one hand, note that we can split the $\dd a$ integral as the sum
\[\int_0^\infty\mr e^{(-\th+\mr iv)a}a^{1-2H_0}\d a
+\int_0^\infty\mr e^{(-\th-\mr iv)a}a^{1-2H_0}\d a.\]
By applying a simple change of variables in the Gamma function identity,
we note that $\int_0^\infty \mr e^{-za}a^{p-1}\d a=\Ga(p)z^{-p}$ whenever
$p>0$ and $\mr{Re}(z)>0$. Given that $\th>0$
and that $1-2H_0=(2-2H_0)-1$, where $2-2H_0>0$ because $H_0<1$,
we conclude that
\[\int_{\mbb R}\mr e^{-\th|a|}\mr e^{\mr iva}|a|^{1-2H_0}\d a
=\Ga(2-2H_0)\big((\th-\mr iv)^{2H_0-2}+(\th+\mr iv)^{2H_0-2}\big).\]
On the other hand, $d-H>0$ and $\th>0$ imply that
\[\int_{\mbb R}\mr e^{-\th|b|}\mr e^{-|v|\cdot|b|}|b|^{d-H-1}\d b
=2\Ga(d-H)(\th+|v|)^{-(d-H)}\]
by the same Gamma function argument.
Thus,
\begin{multline}
\label{Equation: Calculated Integral Formula - H_0<1}
\msf J(\th)
=2\Ga(2-2H_0)\Ga(d-H)\\
\times\int_{-1}^1(1-|v|)\big((\th-\mr iv)^{2H_0-2}+(\th+\mr iv)^{2H_0-2}\big)(\th+|v|)^{-(d-H)}\d v
\end{multline}
when $H_0<1$, and
\begin{align}
\label{Equation: Calculated Integral Formula - H_0=1}
\msf J(\th)
=2\Ga(d-H)\int_{-1}^1(1-|v|)(\th+|v|)^{-(d-H)}\d v
\end{align}
when $H_0=1$.

We now go back to the setting of Lemma \ref{Lemma: Integral Asymptotics - Subcritical},
wherein we assume that $\rho>0$. We aim to apply the dominated convergence theorem in \eqref{Equation: Calculated Integral Formula - H_0<1}
and \eqref{Equation: Calculated Integral Formula - H_0=1}.
Toward this end, we note that whenever $v\neq0$,
\[(\th+|v|)^{-(d-H)}\leq|v|^{-(d-H)},\]
and by the triangle inequality,
\begin{multline*}
\big|(\th-\mr iv)^{2H_0-2}+(\th+\mr iv)^{2H_0-2}\big|\\
\leq|\th-\mr iv|^{2H_0-2}+|\th+\mr iv|^{2H_0-2}
=2|\th^2+v^2|^{H_0-1}
\leq2|v|^{2H_0-2}.
\end{multline*}
Thus, we can justify the application of dominated convergence if we show that
\begin{align}
\label{Equation: Integral Asymptotics - Subcritical 3}
\int_{-1}^1(1-|v|)|v|^{2H_0-2-(d-H)}\d v
=\int_{-1}^1(1-|v|)|v|^{-1+\rho}\d v<\infty.
\end{align}
This holds because $\rho>0$.
We therefore conclude from \eqref{Equation: Calculated Integral Formula - H_0<1} and
\eqref{Equation: Calculated Integral Formula - H_0=1} that
\[\lim_{\th\to0}\msf J(\th)
=2\Ga(2-2H_0)\Ga(d-H)\int_{-1}^1(1-|v|)\big((-\mr iv)^{2H_0-2}+(\mr iv)^{2H_0-2}\big)|v|^{-(d-H)}\d v\]
when $H_0<1$, and that
\[\lim_{\th\to0}\msf J(\th)
=2\Ga(d-H)\int_{-1}^1(1-|v|)|v|^{-(d-H)}\d v\]
when $H_0=1$.
Using the fact that
\[(a\mr i)^b+(-a\mr i)^b=2|a|^b\cos(b\pi/2),
\qquad a,b\in\mbb R,\]
and the simplification
\[\Ga(2-2H_0)\cos\big((H_0-1)\pi\big)=\frac{\Ga(H_0)\Ga(1-H_0)}{2\Ga(2H_0-1)},\qquad 1/2<H_0<1,\]
this simplifies to
\begin{align}
\label{Equation: Integral Asymptotics - Subcritical 4}
\lim_{\th\to0}\msf J(\th)
=\begin{cases}
\displaystyle2\frac{\Ga(H_0)\Ga(1-H_0)}{\Ga(2H_0-1)}\Ga(d-H)\int_{-1}^1(1-|v|)|v|^{-1+\rho}\d v,
&H_0<1,
\vspace{5pt}\\
\displaystyle2\Ga(d-H)\int_{-1}^1(1-|v|)|v|^{-1+\rho}\d v,&H_0=1.
\end{cases}
\end{align}
Given that $\int_{-1}^1(1-|v|)|v|^q\d v=\frac{2}{(2+q)(1+q)}$ whenever $q>-1$, we obtain Lemma \ref{Lemma: Integral Asymptotics - Subcritical}
from \eqref{Equation: Integral Asymptotics - Subcritical 4} by assumption that $\rho>0$.

\subsubsection{Proof of Lemma \ref{Lemma: Integral Asymptotics - Critical}}

We first consider the case $H_0=1$. In that case,
$\rho=0$ implies that $d-H=1$.
Thus, by \eqref{Equation: Calculated Integral Formula - H_0=1}, we get that
\[\msf J(\th)
=4\int_0^1(1-v)(\th+v)^{-1}\d v=4\big((1+\th)\log(1+1/\th)-1\big),\]
where the last equality follows from a direct calculation.

We now consider $1/2<H_0<1$. In that case,
it is challenging to get an asymptotic from \eqref{Equation: Calculated Integral Formula - H_0<1}.
Thus, we consider a different approach.
By applying Fubini's theorem in the first line of \eqref{Equation: J function}, together with the fact that $\rho=0$
(the latter of which implies that $d-H=2H_0-1$),
we can write
\begin{align}
\label{Equation: Integral Asymptotics - Critical TD 2}
\msf J(\th)=
\int_{\mbb R^2}\mr e^{-\th(|a|+|b|)}\left(\int_{-1}^1(1-|v|)\mr e^{\mr iva-|v|\cdot|b|}\d v\right)|a|^{1-2H_0}|b|^{2H_0-2}\d a\dd b.
\end{align}
Given that $z+\bar z=2\mr{Re}(z)$ for any $z\in\mbb C$, we get that
\begin{align}
\label{Equation: -1 to 1 kernel integral 1}
\int_{-1}^1(1-|v|)\mr e^{\mr iva-|v|\cdot|b|}\d v=
2\mr{Re}\left(\int_0^1(1-v)\mr e^{-v(|b|-\mr i a)}\d v\right);
\end{align}
then, an explicit antiderivative calculation yields
\begin{align}
\label{Equation: -1 to 1 kernel integral 2}
2\mr{Re}\left(\int_0^1(1-v)\mr e^{-v(|b|-\mr i a)}\d v\right)=2\mr{Re}\left(\frac{\mr e^{-(|b|-\mr i a)}+(|b|-\mr i a)-1}{(|b|-\mr i a)^2}\right)=:\mathcal K(a,b).
\end{align}
Note that $2\mr{Re}\left(\frac{(|b|-\mr i a)}{(|b|-\mr i a)^2}\right)=\frac{2|b|}{a^2+b^2}$. This is the main singularity of $\mc K(a,b)$. Thus,
we propose to decompose
\begin{align}
\label{Equation: Kernel L}
\mc K(a,b)=\frac{2|b|}{a^2+b^2}+\mc L(a,b),
\qquad\text{where }\mc L(a,b)=2\mr{Re}\left(\frac{\mr e^{-(|b|-\mr i a)}-1}{(|b|-\mr i a)^2}\right),
\end{align}
and, letting $\mc B$ denote the unit ball in $\mathbb R^2$,
use linearity in
\eqref{Equation: Integral Asymptotics - Critical TD 2} to write
\begin{align}
\label{Equation: Integral Asymptotics - Critical TD 3}
\msf J(\th)
=2\int_{\mathcal B^c}\mr e^{-\th(|a|+|b|)}\frac{|a|^{1-2H_0}|b|^{2H_0-1}}{a^2+b^2}\d a\dd b+\msf b(\th),
\end{align}
where
\begin{multline}
\label{Equation: Integral Asymptotics - Critical TD 4}\msf b(\th)=\int_{\mc B^c}\mr e^{-\th(|a|+|b|)}\mc L(a,b)|a|^{1-2H_0}|b|^{2H_0-2}\d a\dd b\\
+\int_{\mc B}\mr e^{-\th(|a|+|b|)}\mc K(a,b)|a|^{1-2H_0}|b|^{2H_0-2}\d a\dd b.
\end{multline}

By a polar change of variables,
\begin{align*}
&2\int_{\mc B^c}\mr e^{-\th(|a|+|b|)}\frac{|a|^{1-2H_0}|b|^{2H_0-1}}{a^2+b^2}\d a\dd b\\
&=2\int_0^{2\pi}|\cos\theta|^{-(2H_0-1)}|\sin\theta|^{2H_0-1}\int_1^\infty\mr e^{-\th r(|\cos\theta|+|\sin\theta|)}r^{-1}\d r\dd\theta\\
&=2\int_0^{2\pi}|\cos\theta|^{-(2H_0-1)}|\sin\theta|^{2H_0-1}\Ga_0\Big(\th\big(|\cos\theta|+|\sin\theta|\big)\Big)\d\theta=4\Ga(H_0)\Ga(1-H_0)\ell(\th),
\end{align*}
where the last line follows by definition of $\ell$ in \eqref{Equation: ell function}. Thus, by \eqref{Equation: Integral Asymptotics - Critical TD 3}, we will
get the first line of \eqref{Equation: Integral Asymptotics - Critical} if we show that the function
$\msf b$ in \eqref{Equation: Integral Asymptotics - Critical TD 4} is continuous and bounded.

For this, by the dominated convergence theorem, it suffices to prove that
\begin{align}
\label{Equation: Integral Asymptotics - Critical TD 5}
&\int_{\mc B}|\mc K(a,b)|\,|a|^{1-2H_0}|b|^{2H_0-2}\d a\dd b<\infty,\\
\label{Equation: Integral Asymptotics - Critical TD 6}
&\int_{\mc B^c}|\mc L(a,b)|\,|a|^{1-2H_0}|b|^{2H_0-2}\d a\dd b<\infty.
\end{align}
Toward this end,
we first claim that there exists
a constant $C>0$ such that
\begin{align}
\label{Equation: -1 to 1 Complex Integral Bound}
|\mc K(a,b)|\leq C\min\big\{1,(a^2+b^2)^{-1/2}\big\}
\end{align}
and
\begin{align}
\label{Equation: L Kernel Upper Bound}
|\mc L(a,b)|\leq C\min\big\{(a^2+b^2)^{-1/2},(a^2+b^2)^{-1}\big\}.
\end{align}
Taking these bounds for
granted for now,
a polar change of variables allows to reduce
\eqref{Equation: Integral Asymptotics - Critical TD 5}
and \eqref{Equation: Integral Asymptotics - Critical TD 6}
to the respective facts that
\[\int_0^1 \min\{1,r^{-1}\}\d r=\int_0^1\d r<\infty
\qquad\text{and}\qquad
\int_1^\infty\min\{r^{-1},r^{-2}\}\d r=\int_1^\infty r^{-2}\d r<\infty.\]
Thus, it only remains to prove
\eqref{Equation: -1 to 1 Complex Integral Bound} and 
\eqref{Equation: L Kernel Upper Bound}.

Regarding \eqref{Equation: -1 to 1 Complex Integral Bound},
the upper bound of 1 is trivial by combining \eqref{Equation: -1 to 1 kernel integral 2} and
\[\left|(1-|v|)\mr e^{\mr iva-|v|\cdot|b|}\right|\leq1,\qquad v\in[-1,1],~a,b\in\mbb R.\]
For the second bound, by the triangle inequality,
\[\left|\frac{\mr e^{-(|b|-\mr i a)}+(|b|-\mr i a)-1}{(|b|-\mr i a)^2}\right|
\leq\frac{1+||b|-\mr i a|+1}{||b|-\mr i a|^2}
=\frac{2+(a^2+b^2)^{1/2}}{a^2+b^2}.\]
In particular,
\[\left|\frac{\mr e^{-(|b|-\mr i a)}+(|b|-\mr i a)-1}{(|b|-\mr i a)^2}\right|
\leq\frac{3(a^2+b^2)^{1/2}}{a^2+b^2}
=\frac{3}{(a^2+b^2)^{1/2}}\qquad\text{when }a^2+b^2\geq1,\]
as desired.
Regarding \eqref{Equation: L Kernel Upper Bound},
by a Taylor expansion,
there exists $C>0$ such that
\[\left|\frac{\mr e^{-(|b|-\mr i a)}-1}{(|b|-\mr i a)^2}\right|
\leq C\left|\frac{(|b|-\mr i a)}{(|b|-\mr i a)^2}\right|\leq\frac{C}{(a^2+b^2)^{1/2}},\qquad \text{when }a^2+b^2\leq 1,\]
and by the triangle inequality,
$\left|\frac{\mr e^{-(|b|-\mr i a)}-1}{(|b|-\mr i a)^2}\right|
\leq\frac{2}{a^2+b^2}.$
With this, the proof of Lemma \ref{Lemma: Integral Asymptotics - Critical} is thus complete.

\subsubsection{Proof of Lemma \ref{Lemma: Integral Asymptotics - Supercritical}}
\label{Section: Proof of Lemma: Integral Asymptotics - Supercritical}

By \eqref{Equation: -1 to 1 kernel integral 1}
and \eqref{Equation: Kernel L},
we can write
\begin{multline*}
\msf J(\th)
=2\int_{\mbb R^2}\mr e^{-\th(|a|+|b|)}\frac{|a|^{1-2H_0}|b|^{d-H}}{a^2+b^2}\d a\dd b\\
+\int_{\mbb R^2}\mr e^{-\th(|a|+|b|)}\mc L(a,b)|a|^{1-2H_0}|b|^{d-H-1}\d a\dd b.
\end{multline*}
This decomposition is convenient since it clearly identifies the power law divergence: By a polar change of variables
and a simple explicit calculation in the radial integral,
\begin{align*}
&2\int_{\mbb R^2}\mr e^{-\th(|a|+|b|)}\frac{|a|^{1-2H_0}|b|^{d-H}}{a^2+b^2}\d a\dd b\\
&=2\int_0^{2\pi}|\cos\theta|^{-(2H_0-1)}|\sin\theta|^{d-H}
\int_0^\infty\mr e^{-\th r(|\cos\theta|+|\sin\theta|)}r^{-1-\rho}\d r\dd\theta\\
&=2\Ga(-\rho)\left(\int_0^{2\pi}\frac{|\cos\theta|^{-(2H_0-1)}|\sin\theta|^{d-H}}{\big(|\cos\theta|+|\sin\theta|\big)^{-\rho}}\d\theta\right)\th^{\rho}\\
&=8\Ga(-\rho)\left(\int_0^{\pi/2}\frac{(\cos\theta)^{-(2H_0-1)}(\sin\theta)^{d-H}}{\big(\cos\theta+\sin\theta\big)^{-\rho}}\d\theta\right)\th^{\rho}.
\end{align*}
Therefore, Lemma \ref{Lemma: Integral Asymptotics - Supercritical} will be
proved if we show that the function
\begin{align}
\label{Equation: z0 function}
\msf z_0(\th)=\int_{\mbb R^2}\mr e^{-\th(|a|+|b|)}\mc L(a,b)|a|^{1-2H_0}|b|^{d-H-1}\d a\dd b
\end{align}
satisfies the claimed properties.

Given that $\th\mapsto \mr e^{-\th(|a|+|b|)}$
is continuous for every $a,b\in\mbb R$, the continuity and boundedness of $\msf z_0$ will follow
from the dominated convergence theorem if we find an integrable function
that dominates the integrand in \eqref{Equation: z0 function} uniformly in $\th$.
Referring back to \eqref{Equation: L Kernel Upper Bound},
we note that
\[\left|\mr e^{-\th(|a|+|b|)}\mc L(a,b)|a|^{1-2H_0}|b|^{d-H-1}\right|
\leq C\min\left\{\frac{|a|^{1-2H_0}|b|^{d-H-1}}{(a^2+b^2)^{1/2}},\frac{|a|^{1-2H_0}|b|^{d-H-1}}{a^2+b^2}\right\}.\]
By a polar change of variables, the integral of the dominating function above is bounded above by
\[C\int_0^1r^{-1-\rho}\d r+C\int_1^\infty r^{-2-\rho}\d r.\]
The first of these integrals is finite when we assume that $\rho<0$.
The second integral is always finite under Assumption \ref{Assumption} because $\rho>-1/2$.

It now only remains to calculate the limit at zero, i.e., \eqref{Equation: Integral Asymptotics - Supercritical Zero Limit}.
By the dominated convergence argument we have just carried out, this is equivalent to
\[\msf z_0(0)=\int_{\mbb R^2}\mc L(a,b)|a|^{1-2H_0}|b|^{d-H-1}\d a\dd b=
\begin{cases}
\displaystyle \frac{4\Ga(H_0)\Ga(1-H_0)\Gamma(d-H)}{\Ga(2H_0-1)(1+\rho)\rho}&H_0>1/2,\\
0&H_0=1/2.
\end{cases}\]
By \eqref{Equation: Kernel L},
we note that the function $\mc L$ is even in both $a$ and $b$.
Indeed, $\mc L(a,-b)=\mc L(a,b)$ is immediate from the fact that $\mc L$ in fact depends only on $|b|$,
and
\[\mc L(-a,b)=2\mr{Re}\left(\frac{\mr e^{-(|b|+\mr i a)}-1}{(|b|+\mr i a)^2}\right)
=2\mr{Re}\left(\frac{\mr e^{-(|b|-\mr i a)}-1}{(|b|-\mr i a)^2}\right)=\mc L(a,b),\]
where the conjugation is absorbed by the real part.
Thus,
\[\msf z_0(0)=4\int_{[0,\infty)^2}\mc L(a,b)a^{1-2H_0}b^{d-H-1}\d a\dd b.\]
Thus, by a polar change of variables,
if we define the function $\msf a(\theta)=\sin\theta-\mr i\cos\theta$ for $\theta\in[0,\pi/2]$, then we get that
\begin{multline*}
\msf z_0(0)=8\mr{Re}\left(\int_0^{\pi/2}\int_0^\infty\frac{\mr e^{-r\msf a(\theta)}-1}{r^2\msf a(\theta)^2}(r\cos\theta)^{1-2H_0}(r\sin\theta)^{d-H-1}r\d r\dd\theta\right)\\
=8\mr{Re}\left(\int_0^{\pi/2}\frac{(\cos\theta)^{1-2H_0}(\sin\theta)^{d-H-1}}{\msf a(\theta)^2}
\left(\int_0^\infty(\mr e^{-r\msf a(\theta)}-1)r^{-\rho-2}\d r\right)\d\theta\right).
\end{multline*}

Given that $-1/2<\rho<0$, the radial part of the integral can be computed by
combining an explicit power-law integral and a Gamma function, whereby
\begin{multline*}
\msf z_0(0)
=8\Gamma(-1-\rho)\mr{Re}\left(\int_0^{\pi/2}\frac{(\cos\theta)^{1-2H_0}(\sin\theta)^{d-H-1}}{\msf a(\theta)^2}
\msf a(\theta)^{1+\rho}\d\theta\right)\\
=8\Gamma(-1-\rho)\mr{Re}\left(\int_0^{\pi/2}(\cos\theta)^{1-2H_0}(\sin\theta)^{d-H-1}\mr e^{\mr i(\rho-1)(\theta-\pi/2)}\d\theta\right),
\end{multline*}
where the second equality follows from $\msf a(\theta)=\sin\theta-\mr i\cos\theta=-\mr i\mr e^{\mr i\theta}=\mr e^{\mr i(\theta-\pi/2)}$
and collecting all the $\msf a(\theta)$ terms. Then, using the change of variables $\phi=\pi/2-\theta$, we get
\[\msf z_0(0)
=8\Gamma(-1-\rho)\mr{Re}\left(\int_0^{\pi/2}(\sin\phi)^{1-2H_0}(\cos\phi)^{d-H-1}\mr e^{\mr i(1-\rho)\phi}\d\phi\right).\]
Our aim is now to transform this integral into a well-known special function, so as to recover the desired representation of $\msf z_0(0)$ as a product of Gamma functions.

Toward this end, we introduce the
change of variables
\[t=\frac{\sin\phi}{\cos\phi+\sin\phi}\in[0,1],\qquad (\sin\phi+\cos\phi)^{2}\dd t=\dd\phi.\]
Given that
\[(\sin\phi)^{1-2H_0}=t^{1-2H_0}(\sin\phi+\cos\phi)^{1-2H_0},\]
\[(\cos\phi)^{d-H-1}=(1-t)^{d-H-1}(\sin\phi+\cos\phi)^{d-H-1},\]
and
\[\mr e^{\mr i(1-\rho)\phi}=
(\cos\phi+\mr i\sin\phi)^{1-\rho}
=(\cos\phi+\sin\phi)^{1-\rho}\big(1-(1-\mr i)t\big)^{1-\rho},\]
we are led to
\[\msf z_0(0)
=8\Gamma(-1-\rho)\mr{Re}\left(\int_0^1\frac{t^{1-2H_0}(1-t)^{d-H-1}}{(1-(1-\mr i)t)^{\rho-1}}(\cos\phi+\sin\phi)^{2(1-\rho)}\d t\right).\]
Finally, noting that
\[t^2+(1-t)^2=\frac{(\sin\phi)^2+(\cos\phi)^2}{(\sin\phi+\cos\phi)^2}=\frac{1}{(\cos\phi+\sin\phi)^2},\]
and moreover, by factoring the roots of that same polynomials, that
\[t^2+(1-t)^2=1-2t+2t^2=\big(1-(1+\mr i)t\big)\big(1-(1-\mr i)t\big),\]
then we obtain a fully $t$-dependent integrand
\[\msf z_0(0)
=8\Gamma(-1-\rho)\mr{Re}\left(\int_0^1\frac{t^{1-2H_0}(1-t)^{d-H-1}}{(1-(1+\mr i)t)^{1-\rho}}\d t\right).\]
By applying \cite[(15.1.2) and (15.6.1)]{NIST:DLMF} and then
\cite[(15.4.6)]{NIST:DLMF}, we get
\begin{multline*}
\msf z_0(0)
=\frac{8\Gamma(-1-\rho)\Ga(2-2H_0)\Ga(d-H)}{\Ga(1-\rho)}\mr{Re}\big({_2}F_1(1-\rho,2-2H_0;1-\rho;1
+\mr i)\big)\\
=\frac{8\Gamma(-1-\rho)\Ga(2-2H_0)\Ga(d-H)}{\Ga(1-\rho)}\mr{Re}\big(\mr i^{2-2H_0}\big)\\
=\frac{8\Gamma(-1-\rho)\Ga(2-2H_0)\Ga(d-H)}{\Ga(1-\rho)}\cos\big(\pi(1-H_0)\big).
\end{multline*}
where ${_2}F_1$ denotes the standard hypergeometric function.
This immediately implies that $\msf z_0(0)=0$ when $H_0=1/2$
since $\cos(\pi/2)=0$. As for $H_0>1/2$,
if we use the trigonometric identity
\[\cos\big(\pi(1-H_0)\big)=\frac{\Ga(H_0)\Ga(1-H_0)}{2\Ga(2-2H_0)\Ga(2H_0-1)},\]
then we obtain the exact formula
\[\msf z_0(0)=\frac{\Ga(H_0)\Ga(1-H_0)}{2\Ga(2-2H_0)\Ga(2H_0-1)}\cdot\frac{8\Gamma(-1-\rho)\Ga(2-2H_0)\Ga(d-H)}{\Ga(1-\rho)}.\]
We then get \eqref{Equation: Integral Asymptotics - Supercritical Zero Limit} through cancellations and
the identity $\Ga(1-\rho)=(1+\rho)\rho\Gamma(-1-\rho)$.

\subsection{Proof of Corollary \ref{Corollary: Limits of E minus Renormalization}}
\label{Section: Limits of E minus Renormalization}

Corollary \ref{Corollary: Limits of E minus Renormalization}-(2) is an immediate consequence of
\eqref{Equation: Renormalization - Supercritical},
the scaling identity in Lemma \ref{Lemma: Scaling of Beta and Alpha},
\eqref{Equation: Separate Expectation Limits - Supercritical}, and
\eqref{Equation: Separate Expectation Limits - Supercritical z limit}.
Thus, we only need to prove Corollary \ref{Corollary: Limits of E minus Renormalization}-(1).

Toward this end,
suppose first that $H_0=1$.
For every $\eps,t>0$, one has
\begin{multline*}
(1+\eps/t)\log(1+t/\eps)
=(1+\eps/t)\log\big((\eps+t)/\eps\big)
=(1+\eps/t)\big(\log(1/\eps)+\log(\eps+t)\big)\\
=(1+\eps/t)\big(\log(1/\eps)-\log\big(1/(\eps+t)\big)\big).
\end{multline*}
Combining this with the Taylor expansion $\log\big(1/(\eps+t)\big)=\log(1/t)+O(\eps)$
and the fact that $(\eps/t)\log(1/\eps)\to0$ as $\eps\to0$, we conclude that
\[(1+\eps/t)\log(1+t/\eps)=\log(1/\eps)-\log(1/t)+o(1)\qquad\text{as }\eps\to0\]
for every fixed $t>0$. Thus, Corollary \ref{Corollary: Limits of E minus Renormalization}-(1)
follows in this case from \eqref{Equation: Renormalization - Critical} and the second line of \eqref{Equation: Separate Expectation Limits - Critical},
the constant $\mf e_\ka$ being the limit of $\ka^2\msf b_2(\th)$ as $\th\to0$.

Suppose then that $H_0<1$.
Well-known asymptotics of the incomplete gamma function (e.g., \cite[5.1.1 and 5.1.11]{AbramowitzStegun})
yields
\[\Ga_0(z)=-\log(z)-\ga_{\mr{EM}}+O(z)\qquad\text{as }z\to0,\]
where $\ga_{\mr{EM}}$ is the Euler-Mascheroni constant.
Given that $|\sin\theta|+|\cos\theta|$ is bounded above and below,
this implies that for every fixed $t>0$,
\[\Ga_0\Big((\eps/t)\big(|\sin\theta|+|\cos\theta|\big)\Big)
=\log(1/\eps)-\log(1/t)-\log\big(|\sin\theta|+|\cos\theta|\big)-\ga_{\mr{EM}}+O(\eps)\]
as $\eps\to0$, where the error is uniform over $\theta\in[0,2\pi]$.
Thus, given that the integral in $\ell(\th)$ is over a bounded set and that
\[\int_0^{2\pi}\left|\frac{\sin\theta}{\cos\theta}\right|^{2H_0-1}\d\theta=2B(H_0,1-H_0)=2\Ga(H_0)\Ga(1-H_0),\]
where $B(\cdot,\cdot)$ denotes the Beta function,
we note that for every $t>0$,
\begin{multline}
\label{Equation: ell function expanded}
\ell(\eps/t)=
\log(1/\eps)-\log(1/t)
-\ga_{\mr{EM}}\\
-\frac{1}{2\Ga(H_0)\Ga(1-H_0)}\int_0^{2\pi}\left|\frac{\sin\theta}{\cos\theta}\right|^{2H_0-1}
\log\big(|\sin\theta|+|\cos\theta|\big)\d\theta
+O(\eps)
\end{multline}
as $\eps\to0$.
Corollary \ref{Corollary: Limits of E minus Renormalization}-(1)
then follows from \eqref{Equation: Renormalization - Critical} and the first line of \eqref{Equation: Separate Expectation Limits - Critical},
the constant $\mf e_\ka$ this time being the sum of the limit of $\ka^2\msf b_1(\th)$ as $\th\to0$ with $-\ka^2\mf b\ga_{\mr{EM}}$
and the integral on the second line of \eqref{Equation: ell function expanded} multiplied by $\ka^2\mf b$.

\subsection{Proof of Corollary \ref{Corollary: SILT=0 when H_0=1/2}}
\label{Section: SILT=0 when H_0=1/2}

Recall the triangular decomposition \eqref{Equation: beta triangle decomposition}.
When $\rho<0$ and $H_0=1/2$, we recall that
the renormalized limit $\ga_t(B)$ is defined as
\eqref{Equation: Renormalized Gamma Triangle Decomposition}. Given that $a_{N,k}^\eps\to a_{N,k}\geq0$
in $L^m$ for all $m\geq1$ as $\eps\to0$
and that $a_{N,1}^\eps,\ldots,a_{N,2^{N-1}}^\eps$
are i.i.d., in order to prove that $\ga_t(B)=0$ it suffices to prove that
\begin{align}
\label{Equation: SILT=0 when H_0=0 1}
\lim_{\eps\to0}\mbf E[a^\eps_{N,1}]=0\qquad\text{for every }N\geq1,
\end{align}
since this will imply that $a_{N,k}=0$ almost surely for all $N$ and $k$.

Toward this end, consider, for every $N\geq1$,
\[\be^\eps_{t/2^{N-1}}(B)=
\int_{[0,t/2^{N-1}]^2}(\msf h\otimes\msf g)_\eps\big(s_2-s_1,B(s_2)-B(s_1)\big)\d s.\]
Decompose the interval above into four squares
\[[0,t/2^N]^2\cup
\big(\underbrace{[0,t/2^N]\times[t/2^N,t/2^{N-1}]}_{J_{N,1}\times I_{N,1}}\big)\cup
\big(\underbrace{[t/2^N,t/2^{N-1}]\times[0,t/2^N]}_{I_{N,1}\times J_{N,1}}\big)\cup
[t/2^N,t/2^{N-1}]^2,\]
where we recall the notations for $J_{N,k}$ and $I_{N,k}$ in \eqref{Equation: J and I Intervals}.
Note that
\[\int_{[0,t/2^N]^2}(\msf h\otimes\msf g)_\eps\big(s_2-s_1,B(s_2)-B(s_1)\big)\d s=\be^\eps_{t/2^N}(B),\]
and
\[\int_{(J_{N,1}\times I_{N,1})\cup(I_{N,1}\times J_{N,1})}(\msf h\otimes\msf g)_\eps\big(s_2-s_1,B(s_2)-B(s_1)\big)\d s
=2a^\eps_{N,1}.\]
Moreover, if we define $\tilde B(s)=B(t/2^{N}+s)-B(t/2^{N})$, which is a standard Brownian motion independent of $B$ restricted to $[0,t/2^N]$,
then we have that
\[\int_{[t/2^N,t/2^{N-1}]^2}(\msf h\otimes\msf g)_\eps\big(s_2-s_1,B(s_2)-B(s_1)\big)\d s=\be^\eps_{t/2^N}(\tilde B).\]
Given that the four squares above intersect in
a region with Lebesgue measure zero, we
get the decomposition
\[\be^\eps_{t/2^{N-1}}(B)
=\be^\eps_{t/2^{N}}(B)
+2a^\eps_{N,1}
+\be^\eps_{t/2^{N}}(\tilde B).\]
Thus, by taking expectations and a trivial rearrangement,
the limit \eqref{Equation: SILT=0 when H_0=0 1} can be reduced to the claim that
\begin{align*}
\lim_{\eps\to0}\big(\mbf E[\be^\eps_{2u}(B)]-2\mbf E[\be^\eps_{u}(B)]\big)=0\qquad\text{for every }u>0
\end{align*}
(in the special case $u=t/2^{N}$).
By Lemma \ref{Lemma: Scaling of Beta and Alpha},
this further reduces to
\begin{align*}
\lim_{\eps\to0}\big((2u)^{1+\rho}\mbf E[\be^{\eps/2u}_{1}(B)]-2u^{1+\rho}\mbf E[\be^{\eps/u}_{1}(B)]\big)=0\qquad\text{for every }u>0.
\end{align*}
By \eqref{Equation: Separate Expectation Limits - Supercritical}, this becomes
\begin{multline*}
\lim_{\eps\to0}\big((2u)^{1+\rho}\big(2\mf d\,(\eps/2u)^{\rho}+\,2\msf z(\eps/2u)\big)-2u^{1+\rho}\big(2\mf d\,(\eps/u)^{\rho}+\,2\msf z(\eps/u)\big)\big)\\
=\lim_{\eps\to0}4u^{1+\rho}\big(2^{\rho}\msf z(\eps/2u)-\msf z(\eps/u)\big)=0,\qquad u>0.
\end{multline*}
This holds by \eqref{Equation: Separate Expectation Limits - Supercritical z limit},
thus concluding the proof of Corollary \ref{Corollary: SILT=0 when H_0=1/2}.

\section{Proof of Proposition  \ref{Proposition: L^2} and Limiting Feynman-Kac Formulas}
\label{Section: L^2}

In this section, our main purpose is to prove
Proposition \ref{Proposition: L^2}. We also take
this opportunity to state the following result,
which is a direct byproduct of that proof:

\begin{proposition}
\label{Proposition: Feynman-Kac}
Let $\ka>0$. Define the function
\begin{align}
\label{Equation: r kappa}
\msf r_\ka(t)=\begin{cases}
0&\rho>0,\\
-\ka^2\mf b\,t\log(1/t)+\mf e_\ka\,t&\rho=0,\\
\tfrac{\ka^2\mf b}{(1+\rho)\rho}t^{1+\rho}&\rho<0,
\end{cases}
\end{align}
where $\mf e_\ka\in\mbb R$ is as in Corollary \ref{Corollary: Limits of E minus Renormalization}.
There exists $\th_\ka>0$ such that for every $t\in(0,\th_\ka)$,
\begin{align}
\label{Equation: Feynman-Kac E}
\mbf E\big[\msf M_{\ka}(t)\big]
&=\mr e^{\msf r_\ka(t)}\int_D\mbf E\left[\mbf 1_{\{\tau_D(B^x)>t\}}\mr e^{\frac{\ka^2}{2}\ga_t(B)}\right]\d x,\\
\label{Equation: Feynman-Kac Var}
\mbf{Var}\big[\msf M_{\ka}(t)\big]
&=\mr e^{2\msf r_\ka(t)}\int_{D^2}\mbf E\Bigg[\mbf 1_{\cap_{i\leq2}\{\tau_D(B^{x_i}_i)>t\}}\mr e^{\sum_{i=1}^2\frac{\ka^2}{2}\ga_t(B_i)}
\left(\mr e^{\ka^2\al_t(x)}-1\right)\Bigg]\d x.
\end{align}
\end{proposition}

We now prove Propositions \ref{Proposition: L^2} and \ref{Proposition: Feynman-Kac}, in that order.

\subsection{Proof of Proposition  \ref{Proposition: L^2}}

It suffices to show that for every $t\in(0,\th_\ka)$,
\begin{align}
\label{Equation: L2 Limit by Completeness 1}
\lim_{\eps_1,\eps_2\to0}\mbf E\big[\msf M_{\ka,\eps_1}(t)\mr e^{-t\ka^2\msf c(\eps_1)}\msf M_{\ka,\eps_2}(t)\mr e^{-t\ka^2\msf c(\eps_2)}\big]
\text{ exists and is finite.}
\end{align}
Consider first the case where $\rho>0$,
wherein $\msf c=0$.
By \eqref{Equation: Finite Epsilon Feynman-Kac 3}, this limit is equal to
\begin{align}
\label{Equation: L2 Limit by Completeness 2}
\lim_{\eps_1,\eps_2\to0}\int_{D^2}\mbf E\left[\mbf 1_{\cap_{i\leq2}\{\tau_D(B_i^{x_i})>t\}}\mr e^{\frac{\ka^2}2\sum_{i=1}^2\be^{\eps_i}_t(B_i)+\ka^2\al_t^{(\eps_1+\eps_2)/2}(x)}\right]\d x.
\end{align}
By a straightforward application of \eqref{Equation: subcritical SILT UI},
\eqref{Equation: MILT UI}, and H\"older's inequality, for every $\ka,t>0$ and $p\geq1$,
one has
\[\sup_{\eps_1,\eps_2>0,~x\in(\mbb R^d)^2}
\mbf E\left[\left(\mbf 1_{\cap_{i\leq2}\{\tau_D(B_i^{x_i})>t\}}\mr e^{\frac{\ka^2}2\sum_{i=1}^2\be^{\eps_i}_t(B_i)+\ka^2\al_t^{(\eps_1+\eps_2)/2}(x)}\right)^p\right]<\infty.\]
Thus, by the dominated convergence and Vitali convergence theorems,
\begin{align*}
\eqref{Equation: L2 Limit by Completeness 2}
=\int_{D^2}\mbf E\left[\mbf 1_{\cap_{i\leq2}\{\tau_D(B_i^{x_i})>t\}}\lim_{\eps_1,\eps_2\to0}\mr e^{\frac{\ka^2}2\sum_{i=1}^2\be^{\eps_i}_t(B_i)+\ka^2\al_t^{(\eps_1+\eps_2)/2}(x)}\right]\d x.
\end{align*}
Then, by \eqref{Equation: subcritical SILT limit} and \eqref{Equation: MILT Limit}, we get that
\begin{align}
\label{Equation: L2 Limit by Completeness 3}
\eqref{Equation: L2 Limit by Completeness 2}
=\int_{D^2}\mbf E\left[\mbf 1_{\cap_{i\leq2}\{\tau_D(B_i^{x_i})>t\}}\mr e^{\frac{\ka^2}2\sum_{i=1}^2\ga_t(B_i)+\ka^2\al_t(x)}\right]\d x,
\end{align}
which concludes the proof of \eqref{Equation: L2 Limit by Completeness 1}.

We now consider the case $\rho\leq0$. If we recall the notation \eqref{Equation: Gamma Epsilon} and
introduce
\begin{align}
\label{Equation: r eps}
\msf r^i_\ka(\eps,t)=
\begin{cases}
0&\rho>0,\\
\big(\mbf E\big[\tfrac{\ka^2}2\be^{\eps}_t(B_i)\big]-t\ka^2\msf c(\eps)\big)&\rho\leq0,
\end{cases}
\end{align}
then by a simple rearrangement of \eqref{Equation: Finite Epsilon Feynman-Kac 3},
the limit in \eqref{Equation: L2 Limit by Completeness 1} is equal to
\begin{align}
\label{Equation: L2 Limit by Completeness 4}
\lim_{\eps_1,\eps_2\to0}\mr e^{\sum_{i=1}^2\msf r^i_\ka(\eps_i,t)}\int_{D^2}\mbf E\left[\mbf 1_{\cap_{i\leq2}\{\tau_D(B_i^{x_i})>t\}}\mr e^{\frac{\ka^2}2\sum_{i=1}^2\ga^{\eps_i}_t(B_i)+\ka^2\al_t^{(\eps_1+\eps_2)/2}(x)}\right]\d x,
\end{align}
On the one hand, by Corollary \ref{Corollary: Limits of E minus Renormalization},
\[\lim_{\eps_1,\eps_2\to0}\mr e^{\sum_{i=1}^2\msf r^i_\ka(\eps_i,t)}=\mr e^{2\msf r_\ka(t)},\]
where $\msf r_\ka(t)$ is defined as in \eqref{Equation: r kappa}. On the other hand,
by a straightforward application of \eqref{Equation: renormalized SILT UI},
\eqref{Equation: MILT UI}, and H\"older's inequality, for every $\ka>0$,
there exists $\th_\ka\in(0,\infty]$ such that
for every $t\in(0,\th_\ka)$, one has
\[\sup_{\eps_1,\eps_2>0,~x\in(\mbb R^d)^2}
\mbf E\left[\left(\mbf 1_{\cap_{i\leq2}\{\tau_D(B_i^{x_i})>t\}}\mr e^{\frac{\ka^2}2\sum_{i=1}^2\ga^{\eps_i}_t(B_i)+\ka^2\al_t^{(\eps_1+\eps_2)/2}(x)}\right)^2\right]<\infty.\]
Thus, by the dominated convergence and Vitali convergence theorems,
\begin{align}
\eqref{Equation: L2 Limit by Completeness 4}
=
\mr e^{2\msf r_\ka(t)}\int_{D^2}\mbf E\left[\mbf 1_{\cap_{i\leq2}\{\tau_D(B_i^{x_i})>t\}}\lim_{\eps_1,\eps_2\to0}\mr e^{\frac{\ka^2}2\sum_{i=1}^2\ga^{\eps_i}_t(B_i)+\ka^2\al_t^{(\eps_1+\eps_2)/2}(x)}\right]\d x
\end{align}
for every $t\in(0,\th_\ka)$.
Then, by \eqref{Equation: renormalized SILT limit} and \eqref{Equation: MILT Limit}, we get that
\begin{align}
\label{Equation: L2 Limit by Completeness 5}
\eqref{Equation: L2 Limit by Completeness 4}
=
\mr e^{2\msf r_\ka(t)}\int_{D^2}\mbf E\left[\mbf 1_{\cap_{i\leq2}\{\tau_D(B_i^{x_i})>t\}}\mr e^{\frac{\ka^2}2\sum_{i=1}^2\ga_t(B_i)+\ka^2\al_t(x)}\right]\d x,
\end{align}
thus concluding the proof of Proposition \ref{Proposition: L^2}.

\subsection{Proof of Proposition \ref{Proposition: Feynman-Kac}}

By Proposition \ref{Proposition: L^2},
\begin{align}
\label{Equation: Feynman-Kac E and Var Limits}
\mbf E\big[\msf M_{\ka}(t)\big]=\lim_{\eps\to0}\mbf E\big[\msf M_{\ka,\eps}(t)\mr e^{-\ka^2t\msf c(\eps)}\big],
\quad
\mbf{Var}\big[\msf M_{\ka}(t)\big]=\lim_{\eps\to0}\mbf{Var}\big[\msf M_{\ka,\eps}(t)\mr e^{-\ka^2t\msf c(\eps)}\big].
\end{align}
With this in hand, the expectation formula \eqref{Equation: Feynman-Kac E} follows from the same argument we used
to calculate the mixed moment limit in \eqref{Equation: L2 Limit by Completeness 1} in the previous section
(whose limit was given in \eqref{Equation: L2 Limit by Completeness 3} and \eqref{Equation: L2 Limit by Completeness 5}), only easier since
there is only one Brownian motion $B$ and the term $\al_t(x)$ is absent. We then get the variance formula
\eqref{Equation: Feynman-Kac Var} by combining the expectation formula
\eqref{Equation: Feynman-Kac E}, the limits for \eqref{Equation: L2 Limit by Completeness 1}
in \eqref{Equation: L2 Limit by Completeness 3} and \eqref{Equation: L2 Limit by Completeness 5},
and the expansion $\mbf{Var}\big[\msf M_{\ka}(t)\big]=\mbf E\big[\msf M_{\ka}(t)^2\big]-\mbf{E}\big[\msf M_{\ka}(t)\big]^2$.

\section{Proof of Theorem \ref{Theorem: Main} Part 1: Expectation Asymptotics}
\label{Section: E}

In this section, we use the tools we have built so far to prove items (1)--(3) in
Theorem \ref{Theorem: Main}. We prove the result case-by-case.

\subsection{Proof of Theorem \ref{Theorem: Main}-(1)}

Suppose that $\rho>0$.
If we combine
\begin{align}
\label{Equation: Exponential Remainder}
\mr e^z=1+z+\msf R(z),\quad z\in\mbb R,
\qquad\text{where }0\leq\msf R(z)\leq z^2(\mr e^z+1),
\end{align}
with the expectation formula \eqref{Equation: Feynman-Kac E}, recalling that $\msf r_\ka(t)=0$ when $\rho>0$, then we get
\begin{multline*}
\mbf E\big[\msf M_{\ka}(t)\big]=\int_D\mbf P\left[\tau_D(B^x)>t\right]\d x\\
+\int_D\mbf E\left[\mbf 1_{\{\tau_D(B^x)>t\}}\tfrac{\ka^2}{2}\ga_t(B)\right]\d x
+\int_D\mbf E\left[\mbf 1_{\{\tau_D(B^x)>t\}}\msf R\big(\tfrac{\ka^2}2\ga_t(B)\big)\right]\d x.
\end{multline*}
Combining this with \eqref{Equation: Heat Content FK}, Theorem \ref{Theorem: Main}-(1) now follows
from two claims: As $t\to0$,
\begin{align}
\label{Equation: Main (1) Dominant}
\int_D\mbf E\left[\mbf 1_{\{\tau_D(B^x)>t\}}\tfrac{\ka^2}{2}\ga_t(B)\right]\d x&\sim\frac{\ka^2\mf b|D|}{(1+\rho)\rho}\,t^{1+\rho},\\
\label{Equation: Main (1) Remainder}
\int_D\mbf E\left[\mbf 1_{\{\tau_D(B^x)>t\}}\msf R\big(\tfrac{\ka^2}2\ga_t(B)\big)\right]\d x&=o(t^{1+\rho}).
\end{align}

For \eqref{Equation: Main (1) Dominant}, given that $\mbf 1_A=1-\mbf 1_{A^c}$, we can write
\[\int_D\mbf E\left[\mbf 1_{\{\tau_D(B^x)>t\}}\tfrac{\ka^2}{2}\ga_t(B)\right]\d x
=|D|\mbf E\big[\tfrac{\ka^2}{2}\ga_t(B)\big]-\int_D\mbf E\left[\mbf 1_{\{\tau_D(B^x)\leq t\}}\tfrac{\ka^2}{2}\ga_t(B)\right]\d x.\]
By combining \eqref{Equation: Gamma Alpha Scalings} and \eqref{Equation: Separate Expectation Limits - Subcritical}, the first term on the right-hand side of this equality has the desired asymptotic. Thus, in order to prove \eqref{Equation: Main (1) Dominant}, it suffices to check that
\begin{align}
\label{Equation: Main (1) Dominant Leftover}
\int_D\mbf E\left[\mbf 1_{\{\tau_D(B^x)\leq t\}}\tfrac{\ka^2}{2}\ga_t(B)\right]\d x=o(t^{1+\rho})\qquad\text{as }t\to0.
\end{align}
By H\"older's inequality and \eqref{Equation: Gamma Alpha Scalings},
\[\text{LHS of }\eqref{Equation: Main (1) Dominant Leftover}
\leq
\tfrac{\ka^2t^{1+\rho}}{2}\mbf E\big[\ga_1(B)^2\big]^{1/2}\int_D\mbf P\left[\tau_D(B^x)\leq t\right]^{1/2}\d x.\]
By \eqref{Equation: subcritical SILT UI}, $\ga_1(B)$'s moments are finite. Thus, to get \eqref{Equation: Main (1) Dominant Leftover},
it suffices to prove that
\begin{align}
\label{Equation: Vanishing of Probability for pleq1}
\text{for every $0<p\leq 1$,}\quad\int_D\mbf P\left[\tau_D(B^x)\leq t\right]^{p}\d x=o(1)\quad\text{as }t\to0
\end{align}
This follows by dominated convergence (see, e.g., \cite[(7.8)]{GaudreauLamarrePan}).

Finally, combining $\mbf 1_A\leq 1$, $0\leq\msf R(z)\leq z^2(\mr e^z+1)$, and \eqref{Equation: Gamma Alpha Scalings},
we get that
\[\text{LHS of }\eqref{Equation: Main (1) Remainder}\leq\frac{\ka^4|D|t^{2(1+\rho)}}{4}\mbf E\left[\ga_1(B)^2\left(\mr e^{\frac{\ka^2}2t^{1+\rho}\ga_1(B)}+1\right)\right].\]
By \eqref{Equation: subcritical SILT UI} and a straightforward application of H\"older's inequality, the expectation in this expression
is bounded above as $t\to0$.
We then obtain \eqref{Equation: Main (1) Remainder} thanks to the fact that $t^{2(1+\rho)}=o(t^{1+\rho})$ as $t\to0$.
With this, the proof of Theorem \ref{Theorem: Main}-(1) is complete.

\subsection{Proof of Theorem \ref{Theorem: Main}-(2) and -(3)}

Consider first the case where $\rho<0$ and $H_0=1/2$.
If we combine \eqref{Equation: Feynman-Kac E} with
the fact that $\msf r_\ka(t)=0$ when $H_0=1/2$ (due to the vanishing
of $\mf b$; see Remark \ref{Remark: b vanishing and rho singularity})
and Corollary \ref{Corollary: SILT=0 when H_0=1/2}, then we get
\[\mbf E\big[\msf M_{\ka}(t)\big]=\int_D\mbf P\left[\tau_D(B^x)>t\right]\d x.\]
This implies that $\mbf E\big[\msf Q_{\ka}(t)\big]=0$ by \eqref{Equation: Heat Content FK}.

Suppose now that $\rho\leq0$ and $H_0>1/2$.
If we apply \eqref{Equation: Exponential Remainder} to
the deterministic exponential
$\mr e^{\msf r_\ka(t)}$ in \eqref{Equation: Feynman-Kac E}, then we get
\begin{multline}
\label{Equation: Main (2)-(3) Decomposition}
\mbf E\big[\msf M_{\ka}(t)\big]=
\int_D\mbf E\left[\mbf 1_{\{\tau_D(B^x)>t\}}\mr e^{\frac{\ka^2}{2}\ga_t(B)}\right]\d x\\
+\Big(\msf r_\ka(t)+\msf R\big(\msf r_\ka(t)\big)\Big)\int_D\mbf E\left[\mbf 1_{\{\tau_D(B^x)>t\}}\mr e^{\frac{\ka^2}{2}\ga_t(B)}\right]\d x.
\end{multline}
By \eqref{Equation: r kappa} and $0\leq\msf R(z)\leq z^2(\mr e^z+1)$, it is clear that
\[\msf r_\ka(t)+\msf R\big(\msf r_\ka(t)\big)\sim\begin{cases}
-\ka^2\mf b\,t\log(1/t)&\rho=0,\\
\frac{\ka^2\mf b}{(1+\rho)\rho}t^{1+\rho}&\rho<0.
\end{cases}\]
Thus, Theorem \ref{Theorem: Main}-(2) and -(3) reduces to
\begin{align}
\label{Equation: Main (2)-(3) Heat Content 1}
\int_D\mbf E\left[\mbf 1_{\{\tau_D(B^x)>t\}}\mr e^{\frac{\ka^2}{2}\ga_t(B)}\right]\d x=\int_Du_0(t,x)\d x+o\big(|\msf r_\ka(t)|\big)\qquad\text{as }t\to0,\\
\label{Equation: Main (2)-(3) Heat Content 2}
\int_D\mbf E\left[\mbf 1_{\{\tau_D(B^x)>t\}}\mr e^{\frac{\ka^2}{2}\ga_t(B)}\right]\d x=|D|+o(1)\qquad\text{as }t\to0;
\end{align}
more specifically, applying \eqref{Equation: Main (2)-(3) Heat Content 1} to the first line of \eqref{Equation: Main (2)-(3) Decomposition},
and then applying \eqref{Equation: Main (2)-(3) Heat Content 2} to the second line of \eqref{Equation: Main (2)-(3) Decomposition}.

By applying \eqref{Equation: Exponential Remainder} to $\mr e^{\frac{\ka^2}{2}\ga_t(B)}$ and using \eqref{Equation: Heat Content FK}, we have that
\[\text{LHS of }\eqref{Equation: Main (2)-(3) Heat Content 1}=\int_Du_0(t,x)\d x
+\int_D\mbf E\left[\mbf 1_{\{\tau_D(B^x)>t\}}\left(\tfrac{\ka^2}{2}\ga_t(B)+\msf R\big(\tfrac{\ka^2}{2}\ga_t(B)\big)\right)\right]\d x.\]
Thus, \eqref{Equation: Main (2)-(3) Heat Content 1} reduces to
\begin{align}
\label{Equation: Main (2)-(3) Heat Content 1 - 1}
\int_D\mbf E\left[\mbf 1_{\{\tau_D(B^x)>t\}}\left(\tfrac{\ka^2}{2}\ga_t(B)+\msf R\big(\tfrac{\ka^2}{2}\ga_t(B)\big)\right)\right]\d x=o\big(|\msf r_\ka(t)|\big)\qquad\text{as }t\to0.
\end{align}
Using once again the fact that $\mbf 1_A=1-\mbf 1_{A^c}$, we get
\begin{multline}
\label{Equation: Main (2)-(3) Heat Content 1 - 2}
\text{LHS of }\eqref{Equation: Main (2)-(3) Heat Content 1 - 1}=
|D|\mbf E\left[\tfrac{\ka^2}{2}\ga_t(B)+\msf R\big(\tfrac{\ka^2}{2}\ga_t(B)\big)\right]\\
-\int_D\mbf E\left[\mbf 1_{\{\tau_D(B^x)\leq t\}}\left(\tfrac{\ka^2}{2}\ga_t(B)+\msf R\big(\tfrac{\ka^2}{2}\ga_t(B)\big)\right)\right]\d x.
\end{multline}
We must now show that both of these expressions are of order $o\big(|\msf r_\ka(t)|\big)$. Given that, when $\rho\leq0$, $\ga_t(B)$ is defined
as the $L^2$ limit of centered random variables (i.e., \eqref{Equation: renormalized SILT limit}), we have that $\mbf E[\ga_t(B)]=0$.
Thus,
\[\left|\mbf E\left[\tfrac{\ka^2}{2}\ga_t(B)+\msf R\big(\tfrac{\ka^2}{2}\ga_t(B)\big)\right]\right|=\mbf E\left[\msf R\big(\tfrac{\ka^2}{2}\ga_t(B)\big)\right]\leq\frac{\ka^4t^{2(1+\rho)}}{4}\mbf E\left[\ga_1(B)^2\left(\mr e^{\frac{\ka^2}2t^{1+\rho}\ga_1(B)}+1\right)\right],\]
where the last inequality follows from $0\leq\msf R(z)\leq z^2(\mr e^{z}+1)$ and \eqref{Equation: Gamma Alpha Scalings}.
This expectation is bounded as $t\to0$ by \eqref{Equation: renormalized SILT UI},
and $t^{2(1+\rho)}=o\big(|\msf r_\ka(t)|\big)$. In particular, the first line of \eqref{Equation: Main (2)-(3) Heat Content 1 - 2} is
of order $o\big(|\msf r_\ka(t)|\big)$. For the second line, a combination of H\"older's inequality,
the scaling \eqref{Equation: Gamma Alpha Scalings}, and $0\leq\msf R(z)\leq z^2(\mr e^{z}+1)$ yields
\begin{multline*}
\left|\int_D\mbf E\left[\mbf 1_{\{\tau_D(B^x)\leq t\}}\left(\tfrac{\ka^2}{2}\ga_t(B)+\msf R\big(\tfrac{\ka^2}{2}\ga_t(B)\big)\right)\right]\d x\right|\\
\leq t^{1+\rho}\mbf E\left[\left(\tfrac{\ka^2}{2}\ga_1(B)+\tfrac{\ka^4}{4}t^{1+\rho}\ga_1(B)^2\left(\mr e^{\frac{\ka^2}2t^{1+\rho}\ga_1(B)}+1\right)\right)^2\right]^{1/2}\int_D\mbf P\left[\tau_D(B^x)\leq t\right]^{1/2}\d x.
\end{multline*}
By \eqref{Equation: renormalized SILT UI}, the expectation above is bounded as $t\to0$.
Combining this with \eqref{Equation: Vanishing of Probability for pleq1} yields
\[\int_D\mbf E\left[\mbf 1_{\{\tau_D(B^x)\leq t\}}\left(\tfrac{\ka^2}{2}\ga_t(B)+\msf R\big(\tfrac{\ka^2}{2}\ga_t(B)\big)\right)\right]\d x=o(t^{1+\rho}).\]
Given that $|\msf r_\ka(t)|\asymp t\log(1/t)$ when $\rho=0$ and $|\msf r_\ka(t)|\asymp t^{1+\rho}$ when $\rho<0$,
we conclude that the second line of \eqref{Equation: Main (2)-(3) Heat Content 1 - 2} is of order $o\big(|\msf r_\ka(t)|\big)$.
This concludes the proof of \eqref{Equation: Main (2)-(3) Heat Content 1}.

With \eqref{Equation: Main (2)-(3) Heat Content 1} established, in order to prove \eqref{Equation: Main (2)-(3) Heat Content 2},
it suffices to show that
\[\int_Du_0(t,x)\d x=|D|+o(1).\]
This follows immediately from \eqref{Equation: Heat Content FK} and \eqref{Equation: Vanishing of Probability for pleq1}
in the case $p=1$. With this, the proof of Theorem \ref{Theorem: Main}-(2) and -(3) is now complete.

\section{Proof of Theorem \ref{Theorem: Main} Part 2: Variance Asymptotic}
\label{Section: Var}

We now conclude the paper with the proof of \eqref{Equation: Variance}.
Note that $\msf M_\ka(t)$ and $\msf Q_\ka(t)$ have the same variance,
so we can in fact prove $\sqrt{\mbf{Var}[\msf M_{\ka}(t)]}\sim \ka\mf c\,t^{H_0}$.

\subsection{Outline}

At first glance, it appears as though
\eqref{Equation: Feynman-Kac Var} is the ideal tool to prove \eqref{Equation: Variance}.
However, for purely technical reasons (e.g., an application of Fubini's theorem
in the Fourier representation \eqref{Equation: Alpha eps Fourier}, which is
easier to justify when the exponential weight coming from the normalization is present),
we work with the prelimit variables. More specifically, recalling the notations for $\ga^\eps_t(B)$ in
\eqref{Equation: Gamma Epsilon} and for $\msf r^i_\ka(\eps,t)$ in \eqref{Equation: r eps},
it follows from the proofs of Propositions \ref{Proposition: L^2} and \ref{Proposition: Feynman-Kac}
that
\begin{align}
\label{Equation: Variance as a Limit}
\mbf{Var}\big[\msf M_\ka(t)\big]
=\lim_{\eps\to0}\mr e^{\sum_{i=1}^2\msf r^i_\ka(\eps,t)}\int_{D^2}\mbf E\left[\mbf 1_{\cap_{i\leq2}\{\tau_D(B^{x_i}_i)>t\}}\mr e^{\frac{\ka^2}2\sum_{i=1}^2\ga^\eps_t(B_i)}\left(\mr e^{\ka^2\al^\eps_t(x)}-1\right)\right]\d x.
\end{align}
We must then provide matching upper and lower bounds for this limit that satisfy the asymptotic \eqref{Equation: Variance}
as $t\to0$.

In informal terms, if $t\approx 0$, then we expect
\[\mr e^{\sum_{i=1}^2\msf r^i_\ka(\eps,t)}\mbf 1_{\cap_{i\leq2}\{\tau_D(B^{x_i}_i)>t\}}\mr e^{\frac{\ka^2}2\sum_{i=1}^2\ga^\eps_t(B_i)}\approx 1,\]
and by a Taylor expansion,
$\mr e^{\ka^2\al^\eps_t(x)}-1\approx\ka^2\al^\eps_t(x).$
Thus, we expect that
\begin{align}
\label{Equation: Variance Proof Heuristic}
\mbf{Var}\big[\msf M_\ka(t)\big]
\approx\lim_{\eps\to0}\ka^2\int_{D^2}\mbf E\left[\al^\eps_t(x)\right]\d x\qquad\text{as }t\to0.
\end{align}
Thanks to this heuristic, we can identify the main mechanism that leads to the variance asymptotics
for $\msf Q_\ka(t)$
as follows:

\begin{lemma}
\label{Lemma: Variance Bound Main Ingredient}
For every bounded open set $G\subset\mbb R^d$,
define
\[\mf c(G)=\left(\int_{G^2}\msf g(x-y)\d x\dd y\right)^{1/2}=\left(\int_{\mbb R^d}|\hat{\mbf 1_G}(\om)|^2\dd\mu(\om)\right)^{1/2};\]
in particular, $\mf c=\mf c(D)$.
For every $m\geq1$,
\[\lim_{\eps\to0}\mbf E\left[\left(\int_{G^2}\al^\eps_t(x)\d x\right)^m\right]^{1/m}\sim \mf c(G)^2\,t^{2H_0}.\]
\end{lemma}

The remainder of this section
is organized as follows: In Section \ref{Section: Lemma: Variance Bound Main Ingredient},
we prove Lemma \ref{Lemma: Variance Bound Main Ingredient}. In Sections \ref{Section: Variance Upper Bound}
and \ref{Section: Variance Lower Bound}, we respectively provide upper and lower bounds
on $\msf M_\ka(t)$'s variance, which, when combined, show a formal version of the approximation \eqref{Equation: Variance Proof Heuristic}
holds (up to smaller order corrections).

\subsection{Proof of Lemma \ref{Lemma: Variance Bound Main Ingredient}}
\label{Section: Lemma: Variance Bound Main Ingredient}

For every $\eps,t>0$, by combining
the Fourier representation \eqref{Equation: Alpha eps Fourier} with the Brownian coupling
\eqref{Equation: Brownian Coupling}, we get
\begin{multline*}
\int_{G^2}\al^\eps_t(x)\d x\\
=\int_{G^2}\int_{[0,t]^2}\int_{\mbb R\times\mbb R^d}\mr e^{-\eps(|\ze|+|\om|^2/2)}\mr e^{\mr i(s_2-s_1)\ze+\mr i\om\cdot(x_2-x_1)+\mr i\om\cdot(B_2(s_2)-B_1(s_1))}\d(\nu\otimes\mu)(\ze,\om)\dd s\dd x.
\end{multline*}
Thanks to the factor $\mr e^{-\eps(|\ze|+|\om|^2/2)}$, we can apply Fubini's theorem to
integrate with respect to $x$ first and then change the order of the remaining integrals, which yields
\begin{multline}
\label{Equation: Variance Bound Main Ingredient 1}
\int_{G^2}\al^\eps_t(x)\d x\\
=\int_{\mbb R^d}\mr e^{-\eps|\om|^2/2}|\hat{\mbf 1_G}(\om)|^2\left(\int_{[0,t]^2}\left(\int_{\mbb R}\mr e^{-\eps|\ze|}\mr e^{\mr i(s_2-s_1)\ze}\d\nu(\ze)\right)\mr e^{\mr i\om\cdot(B_2(s_2)-B_1(s_1))}\d s\right)\dd\mu(\om).
\end{multline}
By \eqref{Equation: mu} and \eqref{Equation: nu}, we can write
\[\int_{\mbb R}\mr e^{-\eps|\ze|}\mr e^{\mr i(s_2-s_1)\ze}\d\nu(\ze)
=\begin{cases}
\frac{\Ga(2H_0+1)}{\Ga(H_0)\Ga(1-H_0)}\int_0^\infty\mr e^{-\eps\ze}\cos\big((s_2-s_1)\ze\big)\ze^{1-2H_0}\d\ze&H_0<1,\\
1&H_0=1,
\end{cases}\]
whenever $s_1\neq s_2$,
where the equality on the right for $H_0<1$ follows from the fact that $\mr e^{\mr i(s_2-s_1)\ze}$'s imaginary part
is odd in $\ze$ (hence the imaginary part of the integral vanishes), and that the corresponding real part is even (hence its
integral over $\mbb R$ is twice the integral over $[0,\infty)$).
If we then apply \cite[3.944-6]{GradshteynRyzhik} (with $\mu=2-2H_0$, $\be=\eps$, and $\de=(s_2-s_1)$),
then we get the exact expression
\begin{align}
\label{Equation: Variance Bound Main Ingredient 2}
\int_{\mbb R}\mr e^{-\eps|\ze|}\mr e^{\mr i(s_2-s_1)\ze}\d\nu(\ze)
=\begin{cases}
\frac{\Ga(2H_0+1)\Ga(2-2H_0)}{\Ga(H_0)\Ga(1-H_0)}\frac{\cos\big((2-2H_0)\arctan\big(\frac{s_2-s_1}{\eps}\big)\big)}{(\eps^2+(s_2-s_1)^2)^{1-H_0}}&H_0<1,\\
1&H_0=1.
\end{cases}
\end{align}
We now analyze the resulting expression on a case-by-case basis.

Suppose first that $H_0=1$. In that case, a combination of \eqref{Equation: Variance Bound Main Ingredient 1}
and \eqref{Equation: Variance Bound Main Ingredient 2} yields
\[\left|\int_{G^2}\al^\eps_t(x)\d x\right|
\leq
\int_{\mbb R^d}\mr e^{-\eps|\om|^2/2}|\hat{\mbf 1_G}(\om)|^2\int_{[0,t]^2}\left|\mr e^{\mr i\om\cdot(B_2(s_2)-B_1(s_1))}\right|\d s\dd\mu(\om)
\leq t^2\mf c(G)^2.\]
Thus, given that $\mr e^{-\eps|\om|^2/2}\to1$ as $\eps\to0$, it follows by dominated convergence that for every $m\geq1$,
one has
\begin{multline}
\label{Equation: Variance Bound Main Ingredient 3}
\lim_{\eps\to0}\mbf E\left[\left(\int_{G^2}\al^\eps_t(x)\d x\right)^m\right]^{1/m}\\
=\mbf E\left[\left(\int_{\mbb R^d}|\hat{\mbf 1_G}(\om)|^2\left(\int_{[0,t]^2}\mr e^{\mr i\om\cdot(B_2(s_2)-B_1(s_1))}\d s\right)\dd\mu(\om)\right)^m\right]^{1/m}.
\end{multline}
By Brownian scaling and the change of variables $u=s/t$, $t^2\dd u=\dd s$, we get
\[\text{RHS of }\eqref{Equation: Variance Bound Main Ingredient 3}=
t^2\mbf E\left[\left(\int_{\mbb R^d}|\hat{\mbf 1_G}(\om)|^2\left(\int_{[0,1]^2}\mr e^{\mr i\om\cdot\sqrt t(B_2(u_2)-B_1(u_1))}\d u\right)\dd\mu(\om)\right)^m\right]^{1/m}.\]
By dominated convergence,
\[\lim_{t\to0}t^{-2}\big(\text{RHS of }\eqref{Equation: Variance Bound Main Ingredient 3}\big)
=\int_{\mbb R^d}|\hat{\mbf 1_G}(\om)|^2\dd\mu(\om)=\mf c(G)^2,\]
thus concluding the proof of Lemma \eqref{Lemma: Variance Bound Main Ingredient} when $H_0=1$.

Now suppose $1/2<H_0<1$.
Given that
\[\left|\frac{\cos\big((2-2H_0)\arctan\big(\frac{s_2-s_1}{\eps}\big)\big)}{(\eps^2+(s_2-s_1)^2)^{1-H_0}}\right|\leq|s_2-s_1|^{2H_0-2}\]
and $\mr e^{-\eps|\om|^2/2}|\mr e^{\mr i\om\cdot(B_2(s_2)-B_1(s_1))}|\leq1$,
By \eqref{Equation: Variance Bound Main Ingredient 1}
and \eqref{Equation: Variance Bound Main Ingredient 2},
there exists a finite constant $C>0$ (which only depends on $H_0$) such that
\[\left|\int_{G^2}\al^\eps_t(x)\d x\right|
\leq C\int_{\mbb R^d}|\hat{\mbf 1_G}(\om)|^2\d\mu(\om)\int_{[0,t]^2}|s_2-s_1|^{2H_0-2}\d s
=\frac{C\mf c(G)^2t^{2H_0}}{H_0(2H_0-1)}<\infty,\]
where the $\dd s$ integral is finite because $H_0>1/2$.
Thus, by dominated convergence,
\begin{multline*}
\lim_{\eps\to0}\mbf E\left[\left(\int_{G^2}\al^\eps_t(x)\d x\right)^m\right]^{1/m}=\frac{\Ga(2H_0+1)\Ga(2-2H_0)\cos\big((2-2H_0)\frac\pi2\big)}{\Ga(H_0)\Ga(1-H_0)}\\
\times\mbf E\Bigg[\Bigg(\int_{\mbb R^d}|\hat{\mbf 1_G}(\om)|^2
\left(\int_{[0,t]^2}|s_2-s_1|^{2H_0-2}\mr e^{\mr i\om\cdot(B_2(s_2)-B_1(s_1))}\d s\right)\dd\mu(\om)\Bigg)^m\Bigg]^{1/m}.
\end{multline*}
By repeating the change of variables/Brownian scaling/dominated convergence argument from the case $H_0=1$
(except that we now also have $|s_2-s_1|^{2H_0-2}$ in the integral,
and we use the calculation $\int_{[0,1]^2}|u_2-u_1|^{2H_0-2}\d u=\frac1{H_0(2H_0-1)}$), we get that
\begin{multline*}
\lim_{\eps\to0}\mbf E\left[\left(\int_{G^2}\al^\eps_t(x)\d x\right)^m\right]^{1/m}\\
\sim\frac{\Ga(2H_0+1)\Ga(2-2H_0)\cos\big((2-2H_0)\frac\pi2\big)}{\Ga(H_0)\Ga(1-H_0)}\cdot t^{2H_0}\mf c(G)^2\cdot\frac1{H_0(2H_0-1)}.
\end{multline*}
We then get the statement of Lemma \ref{Lemma: Variance Bound Main Ingredient} for $1/2<H_0<1$ thanks to the identity
\[\frac{\Ga(2H_0+1)\Ga(2-2H_0)
\cos\big((2-2H_0)\frac{\pi}{2}\big)}
{\Ga(H_0)\Ga(1-H_0)}
=
H_0(2H_0-1).\]

Suppose finally that $H_0=1/2$. In that case, the right-hand side of \eqref{Equation: Variance Bound Main Ingredient 2} reduces to
\[\int_{\mbb R}\mr e^{-\eps|\ze|}\mr e^{\mr i(s_2-s_1)\ze}\d\nu(\ze)=\frac{\eps}{\pi(\eps^2+(s_2-s_1)^2)},\]
since $\dd\nu(\ze)=\dd\mu_{1/2}(\ze)=\tfrac{\dd\ze}{2\pi}$.
Define the function
\[\msf C_\eps(x)=\frac{\eps}{\pi(\eps^2+x^2)},\qquad \eps>0,~x\in\mathbb R,\]
which is the Cauchy density function with scale parameter $\eps$; in particular, one has
$\msf C_\eps(x)=\eps^{-1}\msf C_1(x/\eps)$.
For $i=1,2$, define the map
\[\msf Z_i(\om,s)=\mr e^{\mr i\om\cdot B_i(s)}\mbf 1_{[0,t]}(s),\qquad s\in\mbb R,~\om\in\mbb R^d.\]
Then, we get from \eqref{Equation: Variance Bound Main Ingredient 1}
and \eqref{Equation: Variance Bound Main Ingredient 2} that
\[\int_{G^2}\al^\eps_t(x)\d x
=\int_{\mbb R^d}\mr e^{-\eps|\om|^2/2}|\hat{\mbf 1_G}(\om)|^2\left(\int_{\mbb R}(\msf C_\eps * \msf Z_2(\om,\cdot))(s_1)\overline{\msf Z_1(\om,s_1)}\d {s_1}\right)\dd\mu(\om).\]
Consider an outcome where $B_1$ and $B_2$ are both continuous.
In that case, one has $\msf Z_2(\om,\cdot)\in L^1(\mbb R)$ for every $\om\in\mbb R^d$, and thus
$C_\eps*\msf Z_2\to \msf Z_2(\om,\cdot)$ in $L^1(\mbb R)$
(e.g., \cite[Lemma 2.2.2]{ChenBook}) as $\eps\to0$.
Therefore, given that $\overline{\msf Z_1(\om,\cdot)}\in L^\infty(\mbb R)$ for all $\om\in\mbb R^d$ when $B_1$ is continuous, one has
\[\lim_{\eps\to0}
\int_{\mbb R}(\msf C_\eps * \msf Z_2(\om,\cdot))(s_1)\overline{\msf Z_1(\om,s_1)}\d {s_1}=\int_{\mbb R}\msf Z_2(\om,s)\overline{\msf Z_1(\om,s)}\d s=\int_0^t\mr e^{\mr i\om\cdot (B_2(s)-B_1(s))}\d s\]
for every $\om\in\mbb R^d$. $B_1$ and $B_2$ are continuous almost surely. Thus, almost surely, the above limit
holds pointwise in $\om$.
Moreover, given that $|\msf Z_i(\om,\cdot)|\leq1$, one has
\[\left|\int_{\mbb R}(\msf C_\eps * \msf Z_2(\om,\cdot))(s_1)\overline{\msf Z_1(\om,s_1)}\d {s_1}\right|
\leq\|\msf C_\eps*\msf Z_2(\om,\cdot)\|_{L^1(\mbb R)}
\leq
\|\msf C_\eps\|_{L^1(\mbb R)}\|\msf Z_2(\om,\cdot)\|_{L^1(\mbb R)}=t,\]
where the second inequality follows from Young's convolution inequality.
Therefore, by the dominated convergence theorem,
\[\lim_{\eps\to0}\mbf E\left[\left(\int_{G^2}\al^\eps_t(x)\d x\right)^m\right]^{1/m}
=\mbf E\left[\left(\int_{\mbb R^d}|\hat{\mbf 1_G}(\om)|^2\left(\int_0^t\mr e^{\mr i\om\cdot (B_2(s)-B_1(s))}\d s\right)\dd\mu(\om)\right)^m\right]^{1/m}.\]
By Brownian scaling, and the
change of variables $u=s/t$,
\[\lim_{\eps\to0}\mbf E\left[\left(\int_{G^2}\al^\eps_t(x)\d x\right)^m\right]^{1/m}
=t\mbf E\left[\left(\int_{\mbb R^d}|\hat{\mbf 1_G}(\om)|^2\left(\int_0^1\mr e^{\mr i\om\cdot\sqrt t (B_2(u)-B_1(u))}\d u\right)\dd\mu(\om)\right)^m\right]^{1/m}.\]
We then get the statement of Lemma \ref{Lemma: Variance Bound Main Ingredient} for $H_0=1/2$ by 
dominated convergence. With this, the proof of Lemma \ref{Lemma: Variance Bound Main Ingredient} is now complete.

\subsection{Upper Bound}
\label{Section: Variance Upper Bound}

In this section, we prove that
\begin{align}
\label{Equation: Variance Upper Bound 0}
\mbf{Var}\big[\msf M_\ka(t)\big]\leq
\ka^2\mf c^2\,t^{2H_0}\big(1+o(1)\big)
\qquad\text{as }t\to0.
\end{align}
If we combine the fact that
\begin{align}
\label{Equation: r kappas are negligible}
\lim_{\eps\to0}\mr e^{\sum_{i=1}^2\msf r^i_\ka(\eps,t)}=\mr e^{2\msf r_\ka(t)}=1+o(1)\qquad\text{as }t\to0
\end{align}
with the inequalities
\[\mbf 1_{\cap_{i\leq2}\{\tau_D({_i}B^{x_i})>t\}}\leq1
\qquad\text{and}\qquad
\mr e^z-1\leq z\mr e^z,\quad z\geq0,\]
then we obtain from
\eqref{Equation: Variance as a Limit} that
\begin{align}
\label{Equation: Variance Upper Bound 1}
\mbf{Var}\big[\msf M_\ka(t)\big]
\leq\big(1+o(1)\big)\ka^2\left(\lim_{\eps\to0}\int_{D^2}\mbf E\left[\mr e^{\frac{\ka^2}2\sum_{i=1}^2\ga^\eps_t(B_i)}\al^\eps_t(x)\mr e^{\ka^2\al_t^\eps(x)}\right]\d x\right)
\end{align}
as $t\to0$.
For every $\th>0$, we can split
\begin{multline*}
\int_{D^2}\mbf E\left[\mr e^{\frac{\ka^2}2\sum_{i=1}^2\ga^\eps_t(B_i)}\al^\eps_t(x)\mr e^{\ka^2\al_t^\eps(x)}\right]\d x
=
\int_{D^2}\mbf E\left[\mr e^{\frac{\ka^2}2\sum_{i=1}^2\ga^\eps_t(B_i)}\al^\eps_t(x)\mr e^{\ka^2\al_t^\eps(x)}\mbf 1_{\{\al^\eps_t(x)\leq\th\}}\right]\d x\\
+
\int_{D^2}\mbf E\left[\mr e^{\frac{\ka^2}2\sum_{i=1}^2\ga^\eps_t(B_i)}\al^\eps_t(x)\mr e^{\ka^2\al_t^\eps(x)}\mbf 1_{\{\al^\eps_t(x)>\th\}}\right]\d x.
\end{multline*}
At this point, our aim is to show that for every $\th>0$, one has
\begin{align}
\label{Equation: Variance Upper Bound 2}
\lim_{\eps\to0}\int_{D^2}\mbf E\left[\mr e^{\frac{\ka^2}2\sum_{i=1}^2\ga^\eps_t(B_i)}\al^\eps_t(x)\mr e^{\ka^2\al_t^\eps(x)}\mbf 1_{\{\al^\eps_t(x)\leq\th\}}\right]\d x&\leq\mr e^{\ka^2\th}\mf c^2\,t^{2H_0}\big(1+o(1)\big),\\
\label{Equation: Variance Upper Bound 3}
\lim_{\eps\to0}\int_{D^2}\mbf E\left[\mr e^{\frac{\ka^2}2\sum_{i=1}^2\ga^\eps_t(B_i)}\al^\eps_t(x)\mr e^{\ka^2\al_t^\eps(x)}\mbf 1_{\{\al^\eps_t(x)>\th\}}\right]\d x&=o(t^{2H_0}),
\end{align}
as $t\to0$.
Indeed,
combining this with \eqref{Equation: Variance Upper Bound 1}
and taking $\th\to0$ implies \eqref{Equation: Variance Upper Bound 0}.

By combining
$\mr e^{\ka^2\al_t^\eps(x)}\mbf 1_{\{\al^\eps_t(x)\leq\th\}}\leq\mr e^{\ka^2\th}$
with Tonelli's theorem, H\"older's inequality, and the fact that $\ga^\eps_t(B_1)$ and $\ga^\eps_t(B_2)$ are independent we get
\[\text{LHS of }\eqref{Equation: Variance Upper Bound 2}
\leq\mr e^{\ka^2\th}\lim_{\eps\to0}\mbf E\left[\mr e^{\ka^2\ga^\eps_t(B)}\right]\lim_{\eps\to0}\mbf E\left[\left(\int_{D^2}\al^\eps_t(x)\d x\right)^2\right]^{1/2}.\]
By combining Proposition \ref{Proposition: Limit and UI of SILT} (both the statmements regarding the convergence of $\ga^\eps_t(B)$ and the exponential integrability of $\mr e^{\frac{\ka^2}2\ga^\eps_t(B)}$)
with \eqref{Equation: Gamma Alpha Scalings},
\[\lim_{\eps\to0}\mbf E\left[\mr e^{\ka^2\ga^\eps_t(B)}\right]
=\mbf E\left[\mr e^{\ka^2t^{1+\rho}\ga_1(B)}\right]=1+o(1)\qquad\text{as }t\to0.\]
We thus obtain \eqref{Equation: Variance Upper Bound 2}
by Lemma \ref{Lemma: Variance Bound Main Ingredient}.

Moving on to \eqref{Equation: Variance Upper Bound 3}, by H\"older's inequality,
\begin{multline}\text{LHS of }\eqref{Equation: Variance Upper Bound 3}
\leq|D|^2\sup_{\eps>0}\Bigg(\mbf E\left[\mr e^{2\ka^2\ga^\eps_t(B)}\right]^{1/2}\\
\sup_{x\in(\mbb R^d)^2}\left(\mbf E\big[\al^\eps_t(x)^4\big]^{1/4}
\mbf E\left[\mr e^{4\ka^2\al^\eps_t(x)}\right]^{1/4}\mbf P\big[\al^\eps_t(x)\geq\th\big]^{1/4}\right)\Bigg).
\end{multline}
On the one hand, by \eqref{Equation: subcritical SILT UI},
\eqref{Equation: renormalized SILT UI}, \eqref{Equation: MILT UI}, and Lemma \ref{Lemma: Scaling of Beta and Alpha},
\[\sup_{\eps>0}\Bigg(\mbf E\left[\mr e^{2\ka^2\ga^\eps_t(B)}\right]^{1/2}
\sup_{x\in(\mbb R^d)^2}\left(\mbf E\big[\al^\eps_t(x)^4\big]^{1/4}
\mbf E\left[\mr e^{4\ka^2\al^\eps_t(x)}\right]^{1/4}\right)\Bigg)=O(t^{1+\rho})
\qquad\text{as }t\to0.\]
On the other hand, by Lemma \ref{Lemma: Scaling of Beta and Alpha}, \eqref{Equation: MILT UI}, and Markov's inequality,
for every $p\geq1$,
\[\sup_{\eps>0,~x\in(\mbb R^d)^2}\mbf P\big[\al^\eps_t(x)\geq\th\big]^{1/4}
\leq\sup_{\eps>0,~x\in(\mbb R^d)^2}t^{(1+\rho)p/4}\frac{\mbf E\big[\al^{\eps/t}_1(x)^p\big]^{1/4}}{\th^{p/4}}
=O(t^{(1+\rho)p/4}).\]
Given that $\rho>-1/2$ under Assumption \ref{Assumption},
we can always choose $p$ large enough so that $(1+\rho)(1+p/4)>2H_0$.
Thus, \eqref{Equation: Variance Upper Bound 3} holds, which concludes the proof of \eqref{Equation: Variance Upper Bound 0}.

\subsection{Lower Bound}
\label{Section: Variance Lower Bound}

We now complete the proof of \eqref{Equation: Variance} by showing that
\begin{align}
\label{Equation: Variance Lower Bound 0}
\mbf{Var}\big[\msf M_\ka(t)\big]\geq
\ka^2\mf c^2\,t^{2H_0}\big(1+o(1)\big)
\qquad\text{as }t\to0.
\end{align}
For every $\de>0$, let us denote
\[D_\de=\{x\in D:\mr{dist}(x,\partial D)>\de\}.\]
If we combine \eqref{Equation: r kappas are negligible} with the inequality $\mr e^z-1\geq z$ for $z\geq0$
and a replacement of $D^2$ by $D_\de^2$ in the domain of integration, then we obtain from \eqref{Equation: Variance as a Limit} that
\begin{align}
\label{Equation: Variance Lower Bound 1}
\mbf{Var}\big[\msf M_\ka(t)\big]
\geq\big(1+o(1)\big)\ka^2\left(\lim_{\eps\to0}\int_{D_\de^2}\mbf E\left[\mbf 1_{\cap_{i\leq2}\{\tau_D(B^{x_i}_i)>t\}}\mr e^{\frac{\ka^2}2\sum_{i=1}^2\ga^\eps_t(B_i)}\al^\eps_t(x)\right]\d x\right).
\end{align}

Recall the coupling \eqref{Equation: Brownian Coupling}.
If we define the event
\[A_{t,\de}=\left\{\sup_{0\leq s\leq t}|B_i(s)|<\de\text{ for }i=1,2\right\},\]
then
\[A_{t,\de}\subset\cap_{i\leq 2}\{\tau_D(B^{x_i}_i)>t\}\qquad\text{for every }x_1,x_2\in D_\delta.\]
Thus, by Tonelli's theorem (interchanging the expectation and $\dd x$ integral),
\begin{multline}
\label{Equation: Variance Lower Bound 2}
\lim_{\eps\to0}\int_{D_\de^2}\mbf E\left[\mbf 1_{\cap_{i\leq2}\{\tau_D(B^{x_i}_i)>t\}}\mr e^{\frac{\ka^2}2\sum_{i=1}^2\ga^\eps_t(B_i)}\al^\eps_t(x)\right]\d x\\
\geq\lim_{\eps\to0}\mbf E\left[\mbf 1_{A_{t,\de}}\mr e^{\frac{\ka^2}2\sum_{i=1}^2\ga^\eps_t(B_i)}\int_{D_\de^2}\al^\eps_t(x)\d x\right].
\end{multline}
Finally, by linearity, we can write
\begin{multline}
\label{Equation: Variance Lower Bound 3}
\text{RHS of }\eqref{Equation: Variance Lower Bound 2}=
\lim_{\eps\to0}\mbf E\left[\int_{D_\de^2}\al^\eps_t(x)\d x\right]\\
+\lim_{\eps\to0}\mbf E\left[\left(\mbf 1_{A_{t,\de}}\mr e^{\frac{\ka^2}2\sum_{i=1}^2\ga^\eps_t(B_i)}-1\right)\int_{D_\de^2}\al^\eps_t(x)\d x\right]
\end{multline}

By Lemma \ref{Lemma: Variance Bound Main Ingredient}, for every $\de>0$,
\[\lim_{\eps\to0}\mbf E\left[\int_{D_\de^2}\al^\eps_t(x)\d x\right]\sim \mf c(D_\de)^2\,t^{2H_0}.\]
By the monotone convergence theorem, $\mf c(D_\de)\to\mf c$ as $\de\to0$.
Putting this back into
\eqref{Equation: Variance Lower Bound 2},
\eqref{Equation: Variance Lower Bound 1},
and \eqref{Equation: Variance Lower Bound 3}, it follows that in order to prove
\eqref{Equation: Variance Lower Bound 0},
it suffices to show that the limit on the
second line of \eqref{Equation: Variance Lower Bound 3} is of order $o(t^{2H_0})$ as $t\to0$.
For this, by H\"older's inequality and Lemma \ref{Lemma: Variance Bound Main Ingredient}, it suffices to prove that
\begin{align}
\label{Equation: Variance Lower Bound 4}
\lim_{\eps\to0}\mbf E\left[\left(\mbf 1_{A_{t,\de}}\mr e^{\frac{\ka^2}2\sum_{i=1}^2\ga^\eps_t(B_i)}-1\right)^2\right]^{1/2}=o(1)
\qquad\text{as }t\to0.
\end{align}
By Proposition \ref{Proposition: Limit and UI of SILT} (both the statmements regarding the convergence of $\ga^\eps_t(B)$ and the exponential integrability of $\mr e^{\frac{\ka^2}2\ga^\eps_t(B)}$),
\[\text{LHS of }\eqref{Equation: Variance Lower Bound 4}=\mbf E\left[\left(\mbf 1_{A_{t,\de}}\mr e^{\frac{\ka^2}2\sum_{i=1}^2\ga_t(B_i)}-1\right)^2\right]^{1/2}\]
for small enough $t>0$.
For every $\de>0$,
it follows from Brownian scaling (including its application to $\ga_t(B)$ in \eqref{Equation: Gamma Alpha Scalings}) that
\[\lim_{t\to0}\mbf 1_{A_{t,\de}}\mr e^{\frac{\ka^2}2\sum_{i=1}^2\ga_t(B_i)}=1\qquad\text{in probability}.\]
Therefore, 
\eqref{Equation: Variance Lower Bound 4} holds thanks to the uniform integrability
estimates \eqref{Equation: subcritical SILT UI} and \eqref{Equation: renormalized SILT UI}. With this, the proof
of \eqref{Equation: Variance Lower Bound 0},
and therefore of
\eqref{Equation: Variance}, is complete.

\bibliographystyle{plain}
\bibliography{Bibliography}

\end{document}